\documentclass[12pt]{article}
\usepackage[a4paper,margin=2cm]{geometry}

\usepackage{amsfonts,graphicx,amsmath,amssymb,amsthm,color,nicefrac,enumerate, layout, bbm, listings,  hyperref}
\usepackage[latin1]{inputenc}
\usepackage{natbib}
\numberwithin{equation}{section}

\newtheorem{lemma}{Lemma}
\newtheorem{theorem}{Theorem}

\theoremstyle{definition}

\theoremstyle{remark}

\newtheorem{remark}{Remark}

\newcommand{\convas}[1][]{\xrightarrow[#1]{\mathrm{a.s.}}}
\newcommand{\convp}[1][]{\xrightarrow[#1]{\PP}}
\newcommand{\convl}[1][]{\xrightarrow[#1]{\mathrm{law}}}

\newcommand{\PP}{\mathbf{P}}
\newcommand{\EE}{\mathbf{E}}
\newcommand{\QR}{\mathbf{Q}}

\newcommand{\NN}{\mathbb{N}}
\newcommand{\RR}{\mathbb{R}}
\newcommand{\X}{\mathbb{X}}

\newcommand{\cL}{\mathcal{L}}

\newcommand{\cN}{\mathcal{N}}

\newcommand{\cM}{M}

\newcommand{\cQ}{Q}

\newcommand{\ind}[1]{\mathbf{1}_{#1}}

\newcommand{\Linfty}{L_\infty}

\newcommand{\vd}{\,\mathrm{d}}
\newcommand{\dd}{\mathrm{d}}

\DeclareMathOperator{\sgn}{sgn}

\def\eqlaw{\stackrel{\mathrm{law}}{=}}
\newcommand{\eqdef}{\mathbin{:=}}

{\begin{itemize}[leftmargin=1cm,noitemsep,label={$\star$}]}{\end{itemize}}

\begin{document}

\title{
Estimation of parameters and local times in a discretely observed threshold diffusion model
}
\author{
Sara Mazzonetto\footnote{Universit\'e de Lorraine, CNRS, Inria, IECL, F-54000 Nancy, France. E-mail: \texttt{sara.mazzonetto@univ-lorraine.fr}}
\ and
Paolo Pigato\footnote{%
Department of Economics and Finance, University of Rome Tor Vergata, Via Columbia 2, 00133 Roma, Italy. 
E-mail: \texttt{paolo.pigato@uniroma2.it
} 
}
}

\date{\today}

\maketitle

\begin{abstract}
\noindent
We consider a simple mean reverting diffusion process, with piecewise constant drift and diffusion coefficients,  discontinuous at a fixed threshold. We discuss estimation of drift and diffusion parameters from discrete observations of the process, with a generalized moment estimator and a maximum likelihood estimator. We develop the asymptotic theory of the estimators when the time horizon of the observations goes to infinity, 
considering both cases of a fixed time lag (low frequency) and a vanishing time lag (high frequency) between consecutive observations. In the setting of low frequency observations and infinite time horizon we also study the convergence of three local time estimators, that are already known to converge to the local time in the setting of high frequency observations and fixed time horizon. We find that these estimators can behave differently, depending on the assumptions on the time lag between observations. 
\end{abstract}

\bigskip


\medskip

\noindent{\textbf{AMS 2020: }} 
	primary: 
		62F12; 
	secondary: 	
		62M05; 
		60F05; 
		60J55. 


\section{Introduction}

We consider the diffusion process solution to the following stochastic differential equation (SDE)
\begin{equation}
    \label{eq:DOBM}
    X_t=X_0+\int_0^t \sigma(X_s)\vd W_s+\int_0^t b(X_s) \vd s , \quad t\geq 0,
\end{equation}
with piecewise constant volatility and drift coefficient, possibly discontinuous at $r\in \RR$,
\begin{equation*}
    \label{sigmaDOBM}
    \sigma(x)=\begin{cases}
	\sigma_+>0		&\text{ if }x\geq r,\\
	\sigma_->0		&\text{ if }x<r,
    \end{cases}
\quad \text{and} \quad 
    b(x)=\begin{cases}
	b_+ \in \RR 	&\text{ if }x\geq r,\\
	b_-  \in \RR	&\text{ if }x<r.
    \end{cases}
\end{equation*}
We assume the initial condition to be deterministic. 
Separately on $(r,\infty)$ and $(-\infty,r)$, the process follows the drifted Brownian motion dynamics with different parameters.
We also assume that $b_- >0$ and $b_+ < 0$, so that the process is mean reverting and ergodic.

We discuss the estimation of parameters from discrete observations of $X$, when the threshold $r$ is known. We consider two types of drift estimator: a generalized moment estimator (GME) and a maximum discretized likelihood estimator (dMLE). We prove that they are equivalent in a suitable sense, so that  the asymptotic theory we develop applies to both.

We construct the GME, based on discrete observations, using the theoretical long time behavior of conditional moments and local times of $X$.
Maximum likelihood estimation (MLE) of this process has been considered in \cite{lp2}, assuming that $X$ is continuously monitored and the observation span goes to infinity. Since we assume here to observe the process at discrete times, we study the corresponding dMLE, with observation span going to infinity.

We consider two types of asymptotic setting: (a) fixed time lag $h$ (low frequency observations) and number of observations $N$ going to infinity, and (b) shrinking time lag $h_{N}$ (high frequency observations) and simultaneously time horizon $T_{N}=Nh_{N}$ going to infinity. In all cases we prove consistency and asymptotic normality of GME and dMLE estimators. We also show that the dMLE based on $N$ observations converges in high frequency to the MLE based on continuous observations, with speed $N^{1/4}$, and an analogous result for the GME.
Moreover, in setting (a), we propose a GME for the diffusion coefficient, for which we prove again consistency and asymptotic normality. 

Finally, we consider three different estimators for the local time at the discontinuity level $r$, that are known to converge to the local time with speed $N^{{1/4}} $ in the high frequency limit of the observations. 
We show that when the observation lag is fixed and the time horizon goes to infinity, the asymptotic behavior is the same as for the local time only for one of the estimators, while the other two have different limit behaviors, that depend on the observation lag.

\paragraph{Related work.}
The solution to \eqref{eq:DOBM} with piecewise affine instead of piecewise constant drift is the so-called threshold Ornstein-Uhlenbeck (OU) process, for which parameters inference (including the threshold $r$) has been discussed e.g. in~\cite{kutoyants2012, dieker2013,  su2015,su2017, YU2020,  HuXi, mazzonetto2020drift}. Estimation of a piecewise constant drift as in \eqref{eq:DOBM}, despite the process having several applications (see below), has been considered less, see e.g. \cite{motaesquivel, lp2}.

\smallskip

In the low frequency observations setting, the main tool we use in this paper is the ergodic theorem in the form of \cite[Theorem 17.0.1]{MeynTweedie}, in the multidimensional framework as in \cite[Section 1.8]{brooks2011handbook}, in the spirit of~\cite{HuXi}. 
In the high frequency observations setting, we obtain that the asymptotic behavior is the same as for the continuous observations case, under the condition $\lim_{N\to\infty} h_N T_N =0$, analogous to the condition given for diffusions with no threshold in \citep{kessler,alaya_kebaier_2013,amorino_glotier} and with threshold in~\cite{mazzonetto2020drift}.

\smallskip

As already mentioned, several local time estimators for Brownian motion or skew Brownian motion, with fixed time horizon and $N$ observations, have been considered in the literature (see e.g.~\cite{port1,Borodin,LP,lp2,lmt1}) and they are known to converge to the local time with speed $N^{{1/4}} $ in the high-frequency limit of the observations \citep{j1,mazzonetto2019rates}. 
Local time approximation is an important topic in statistics of processes, for instance because the amount of time spent at a certain level is related to the accuracy of the estimation of the process at that level \citep{florens-zmirou}.
Recently in \cite{ChristensenStrauchAAP,ChristensenStrauchTrottnerBernoulli} the local time of a scalar diffusion model have been used in the exploration phase of studies of the exploration vs exploitation tradeoff in a reinforcement learning setting. The authors point out the difficulties related to local time estimation and propose strategies to avoid it.

\smallskip

The process in \eqref{eq:DOBM} and variations of it have been widely used in financial modelling, from the point of view of time series (see e.g. \cite{ang,motaesquivel,lp1}), options pricing and implied volatility (see e.g. \cite{liptonsepp,gairat,DongWong,lipton:2018,pigato,buckner}), interest rates (see e.g. \cite{pai,interestrate,su2015,su2017}) and others.
The discrete time analog of threshold diffusions are Threshold autoregressive (TAR) and in particular self-exciting (SETAR) models \citep{Tong:2011ud,Chen:2011bk}. They have been widely used in financial modelling, recently also in combination with reinforcement learning in \cite{ghp2023}.
Threshold processes have also a wide range of applications outside of financial modelling; for example,
\cite{Hottovy}  discusses deterministic and stochastic triggers in threshold diffusion models for rainfall and convection. For other applications of the specific threshold diffusion in \eqref{eq:DOBM} we refer to \cite{lp2}.

\paragraph{Outline.} 
We state our main results in Section~\ref{sec:results} and we discuss them in Section~\ref{rm:derivation}.
We present some numerical studies  in Section~\ref{sec:numerics}.
We collect proofs and technical material in Section~\ref{sec:proofs}.

\section{Estimators and convergence results} \label{sec:results}

Without loss of generality, we assume the threshold is at $r=0$. Indeed, given $(X_t)_{t}$ solution to~\eqref{eq:DOBM} with known level (threshold) $r\in \RR$, the process $(Y_t)_t:=(X_t-r)_{t}$ solves~\eqref{eq:DOBM} with threshold at $0$ and initial condition $Y_0=X_0-r$.
As a consequence, all the results in this document can be easily rewritten in the case $r\neq 0$.

So, let $X=(X_{t})_{t\in [0,\infty)}$ be the solution to~\eqref{eq:DOBM} with $r=0$, $W$ Brownian motion, and $X_0$ deterministic (see~\cite{legall} or~\cite{bass2005one} for strong existence results). We assume that $b_- >0$ and $b_{+}<0$, so that the process is mean reverting and ergodic.
In this case, also the discrete time process $(X_{hk})_{k\in\NN}$ (the so-called $h-$skeleton sampled chain) is ergodic, for any $h>0$ (this can be proved as in \cite[Lemma 1]{HuXi}). 

Given $T\in (0,\infty)$, we assume to observe $X$ over a time grid with $h$ discretization step such that $T=hN$.
The (symmetric) local time of $X$ at $0$, denoted by $L_T(X)$, is a continuous process quantifying the amount of time spent by $X$ close to level $0$ up to time $T$. An estimator of the local time $L_T(X)$ is given by
\begin{equation} \label{def:LT:discr}
	L_{h,N}:= 2
	 \sum_{k=0}^{N-1} \ind{\{X_{kh
	 } X_{(k+1)h}<0\}}|X_{(k+1)h}|.
\end{equation}
We refer to Section~\ref{sec:local_time} for details on the local time, convergence results and alternative local time estimators.
Writing $\pm\in \{-,+\}$, let us introduce quantities
\begin{equation}\label{def:M:Q}
 \cM_T^{\pm}:=
	 \int_0^T  \ind{\{ \pm X_s \geq 0 \} } \vd X_s
\quad \text{and} \quad	
\cQ_T^{\pm} 	
:=
\int_0^T \ind{\{ \pm X_s \geq 0 \}}\vd s
\end{equation}
and their discrete counterparts
\begin{equation} \label{def:M:Q:disc}
    \cM^{\pm}_{h,N} := \sum_{k=0}^{N-1} 
    \ind{ \{ \pm X_{kh} \geq 0 \} } (X_{(k+1)h}-X_{kh})
\quad \text{and} \quad
    \cQ^{\pm}_{h,N} := h \sum_{k=0}^{N-1} 
    \ind{\{ \pm X_{kh} \geq 0 \}}.
\end{equation}
Throughout the paper we use the notation $\cN(\mu,\Sigma)$ for a Gaussian variable with mean vector $\mu$ and covariance matrix $\Sigma$. We also write $\cN=\cN(0,\operatorname{Id})$ for a standard Gaussian. We use in what follows the notion of \emph{stable convergence}, that we denote $\xrightarrow{\mathrm{stably}}$, for which we refer to \cite{renyi,js,jp}. We also write $\convas$ for almost sure convergence, $\convp$ for convergence in probability, $\convl$ for convergence in law.

\subsection{Drift parameters inference}

We consider estimators for the drift parameters $b_+,b_-$ in \eqref{eq:DOBM}, based on the quantities above (see next Section \ref{rm:derivation} for a derivation). The first one is a GME
\begin{equation}\label{eq:GME}
\begin{pmatrix}
	 \overline{b}^+_{h,N}, & \overline{b}^-_{h,N}
	\end{pmatrix}
	= 
	\begin{pmatrix}
- \frac{L_{h,N}}{2\cQ_{h,N}^+}, & 
\frac{L_{h,N}}{2\cQ_{h,N}^-}	\end{pmatrix}.
\end{equation}
The second one is a dMLE, that can also be interpreted as a least squares estimator (LSE)
\begin{equation}\label{eq:QMLE_dis_time}
	\begin{pmatrix}
	 \widehat{b}^+_{h,N}, &  \widehat{b}^-_{h,N}
	\end{pmatrix}
	= 
	\begin{pmatrix}
	\frac{	\cM^{+}_{h,N} }{ 
\cQ^{+}_{h,N} }
, & \frac{	\cM^{-}_{h,N} }{ 
\cQ^{-}_{h,N} } 
	\end{pmatrix}.
	\end{equation}
In what follows, we provide a complete asymptotic theory for these estimators.
We consider different settings, depending on the assumptions on $h$ and $T$, which can be as follows:
\begin{enumerate}
\item $h$ is a fixed constant and $N\to \infty$ (\emph{long time}, since $T=N h$, and \emph{low frequency});
\item $h=h_{N}=T/N$, with $T$ fixed and $N\to \infty$ (\emph{high frequency});
\item $h=h_{N}=T_N/N$, with $T_N \to \infty$ and $h_N \to 0$ as $N\to\infty$ (\emph{long time} and \emph{high frequency}). 
\end{enumerate}
We show that the estimators in \eqref{eq:GME} and \eqref{eq:QMLE_dis_time} are \emph{asymptotically equivalent} as the time horizon goes to infinity (see Lemma~\ref{lemma:equivalence} in Section~\ref{sec:relation}), so that both next Theorems~\ref{th:fixed:h} and~\ref{th:joint:CLT} hold for both estimators, in exactly the same form. For this reason, we state our results in extended form only for $\overline{b}^\pm$. 

We first consider the long time, low frequency framework, meaning that we observe the process over a discrete time grid with fixed time lag, with number of observations going to infinity.

\begin{theorem} \label{th:fixed:h}
Let $h>0$ fixed. Then
\begin{enumerate}[i)]
\item the estimator is strongly consistent:
$
   \ \overline{b}^\pm_{h,N}
\convas[N\to\infty]
   \ b_\pm
$,
\item
the estimator is asymptotically normal, i.e.
\[
\sqrt{N}
\begin{pmatrix}
 \overline{b}^+_{h,N}
-b_+ \\
 \overline{b}^-_{h,N}
-b_-
\end{pmatrix}
\xrightarrow[N\to \infty]{\mathrm{law}}
 \cN(0, \Gamma),
\]
where $\Gamma$ is given in \eqref{eq:def:gammaq} in Lemma~\ref{lem:jointconv}.
\end{enumerate}
These results hold also substituting estimator $\overline{b}^\pm$ with $\widehat{b}^\pm$, but in this case the convergence in $i)$ holds in probability (weak consistency).
\end{theorem}
The result for $\overline{b}^\pm$ is proved in Section~\ref{sec:gme} using ergodic theorems.
By Lemma \ref{lemma:equivalence} in Section~\ref{sec:relation}, we deduce the result for $\widehat{b}^\pm$ from the one for $\overline{b}^\pm$. 

\medskip

We now assume to observe the process over a discrete time grid with vanishing time lag, with the number of observations going to infinity in a way so that the time horizon also goes to infinity (long time and high frequency framework).

\begin{theorem} \label{th:joint:CLT}
Let $(T_N)_{N\in \NN},(h_N)_{N\in \NN}$ be positive sequences satisfying $h_N=T_N/N$ and
\[
\lim_{N\to\infty}T_N = \infty \quad \text{and} \quad \lim_{N\to\infty} h_N = 0.
\]
Then
\begin{enumerate}[i)]
\item the estimator is consistent: 
$
   \ \overline{b}^\pm_{h_N,N}
\convp[N\to\infty]
   \ b_\pm,
$
\item \label{th2:item:CLT}
if in addition 
\begin{equation} \label{cond:TN}
 \lim_{N\to\infty} T_N h_N =0,
\end{equation}
the estimator is asymptotically normal, i.e.
\[
\sqrt{T_N}
\begin{pmatrix}
 \overline{b}^+_{h_N,N}
-b_+ \\
 \overline{b}^-_{h_N,N}
-b_-
\end{pmatrix}
\xrightarrow[N\to \infty]{\mathrm{stably}} \cN(0, \Sigma)
\]
where
\[ \Sigma  
	:= \frac{b_-+|b_+|}{b_-|b_+|}\begin{pmatrix} \sigma_+^{2} |b_+| & 0 \\
		0 & \sigma_-^{2} b_- \end{pmatrix}.
\] 
\end{enumerate}
Precisely the same results hold substituting estimator $\overline{b}^\pm$ with $\widehat{b}^\pm$.
\end{theorem}
\medskip

The result for $\widehat{b}^\pm$ is proved in Section \ref{sec:mle}. The result for $\overline{b}^\pm$ follows from Lemma \ref{lemma:equivalence} in Section~\ref{sec:relation}.

The next result shows that when the time horizon is fixed, both the estimators above, based on discrete observations, converge  in high frequency, with speed $N^{1/4}$, to their continuous time analogues.
\begin{theorem}\label{th:disc}
    Let $T\in (0,\infty)$ be fixed and $h_N=T/N$. 
Then, 
both
\begin{equation}
N^{1/4}
	(\widehat{b}^+_{h_N,N} -  
	\frac{	\cM^{+}_T }{ Q^{+}_T }
	, \widehat{b}^-_{h_N,N} - 
	\frac{	\cM^{-}_T }{ Q^{-}_T } )  
\qquad \text{ and } \qquad 
N^{1/4}
	\left(\overline{b}^+_{h_N,N} -  \frac{-L_T}{2Q^{+}_T} , \overline{b}^-_{h_N,N} - 	\frac{L_T}{ 2Q^{-}_T } \right),  
\end{equation}
when $N\to\infty$, converge stably to 
\begin{equation} 
	\left( 	
	  \frac1{Q^{+}_{T}}, -  \frac1{Q^{-}_{T}}	
 \right) 
	 \sqrt{\frac{4 \sqrt{T}}{3 \sqrt{2 \pi}}  \frac{\sigma_-^2 +\sigma_+^2}{\sigma_-+\sigma_+}} {\sqrt{L_T(X)} } \  \cN 
\end{equation}
with $\cN$ a standard Gaussian random variable independent of $X$, and $L_T(X)$  local time of $X$ at $0$, up to time $T$.
\end{theorem}

\medskip

This theorem is proved in Section \ref{pr:thdisc}.

\subsection{Volatility parameters inference}
\label{sec:volat:estim}
We now introduce an estimator for the volatility parameters $\sigma_\pm$ based on discrete observations with a fixed $h$ discretization step. Let us denote 
\begin{equation}\label{def:Q1:disc}
Q^{\pm,1}_{h,N} := h \sum_{k=0}^{N-1}  X_{kh}
    \ind{\{ \pm X_{kh} \geq 0\}}
\end{equation}
and
\begin{equation}\label{eq:def:est:sigma}
(\bar{\sigma}^\pm_{h,N})^{2}
=
\pm
\frac12 \frac{L_{h,N} Q^{\pm,1}_{h,N}}{ ({Q}^{\pm}_{h,N})^2} .
\end{equation}
This is a GME estimator as the one in \eqref{eq:GME}. We consider its long time behavior. 

\begin{theorem} \label{th:sigma:h}
Let $h>0$ fixed. Then
\begin{enumerate}[i)]
\item the estimator is consistent:
$
   \  \bar{\sigma}^\pm_{h,N}
\convp[N\to\infty]
   \sigma_\pm,
$
\item
the estimator is asymptotically normal, i.e.
\[
\sqrt{N}
\begin{pmatrix}
( \bar{\sigma}^+_{h,N})^{2}
-\sigma^{2}_+ \\
 ( \bar{\sigma}^-_{h,N})^{2}
-\sigma^{2}_-
\end{pmatrix}
\xrightarrow[N\to \infty]{\mathrm{law}}
  \cN(0, \Gamma_{\sigma})
\]
where $\Gamma_{\sigma}$ is given in \eqref{eq:def:gammav}.
\end{enumerate}
\end{theorem}
This theorem is proved in Section \ref{sec:gme}. 

When high frequency observations are available, the quadratic variation estimators considered in \cite{LP} seem more suitable for estimating $\sigma_{\pm}$. However, $\bar{\sigma}^\pm_{h,N}$ in \eqref{eq:def:est:sigma} could be the better choice for estimating $\sigma_{\pm}$ when $X$ is observed at a low frequency, since high frequency observations are not required for convergence. Let us also note that a joint version of Theorems \ref{th:fixed:h} and \ref{th:sigma:h} could be written, providing a joint convergence of $\bar{b}^\pm_{h,N}, \bar{\sigma}^\pm_{h,N}$, with speed $\sqrt{N}$ and asymptotic covariance again of the form in \eqref{eq:def:gamma2}, this time a $4$ by $4$ matrix.

\subsection{Local time estimation}
\label{sec:local_time}
Given a semi-martingale $Y \colon \Omega \times [0,\infty) \to \RR$
and $z\in \RR$, $T\in [0,\infty)$, the quantity 
\begin{equation*}
    L_T^z(Y)=\lim_{\epsilon\to 0}\frac{1}{2\epsilon}\int_0^T
    \ind{\{-\epsilon\leq Y_s-z\leq \epsilon\}}d\langle Y\rangle_s
\end{equation*}
defines the (symmetric) local time of $Y$ at $z$, which quantifies the amount of time spent by $Y$ close to level $z$ up to time $T$, properly re-scaled. 
We focus here on the local time at the discontinuity of $X$, that is $L_{T}:=L_T^0(X)$. Because of the ergodic property, the rescaled local time converges: 
\begin{equation}\label{eq:lt:limit}
\frac{L_{T}(X)}{T} \convas[T\to \infty]  \frac{2 |b_+| b_-}{b_-+|b_+|} =:\Linfty 
\end{equation}
(cf.~\cite[Chapter II 35(c)]{borodin2015handbook}).
There are several ways to approximate the local time $L_T$ from discrete observations of $X$. 
We have already introduced the estimator ${L}_{h,N}$ in \eqref{def:LT:discr}. 
A well known estimator is the one from renormalization of number of crossings, which is 
\[
\bar{L}_{h,N}=	\sqrt{\frac{
\pi}{2}}\frac{\sigma_++\sigma_-}{2}
  \sqrt{h}\sum_{k=0}^{N-1} \ind{\{X_{kh
	 } X_{(k+1)h}<0\}}.
\] 
We define also the estimator ``from sample covariance''
\[
\hat{L}_{h,N}=	
-\frac{3\sqrt{\pi}}{2\sqrt{2}}
\frac{\sigma_-+\sigma_+}{\sigma_-\sigma_+}
  \frac{1}{\sqrt{h}}\sum_{k=0}^{N-1} [X^{+},X^{-}]^{h}_{N}
  \]
where $[Y,Z]^{h}_{N}=\sum_{k=0}^{N-1} (Y_{(k+1)h}-Y_{kh}) (Z_{(k+1)h}-Z_{kh})$ is the discrete covariation of two one-dimensional processes $Y,Z$, over $N$ observations with a time lag $h$.
The process $X^{+}=X \ind{\{X \geq 0 \}}$ denotes the positive part of $X$ and $X^{-}=- X \ind{\{X< 0 \}}$ its negative part (which is positive). 

The three statistics above converge to the local time for fixed final horizon $T=h_{N} N$ when $N$ (the number of observations) goes to infinity, with speed $N^{1/4}$ (see~\cite{mazzonetto2019rates}).
We now study their behavior for low frequency (fixed time lag) observations, in long time.

Let us denote with $O(h)$ any function $f(h)$ s.t.~$f(h)/h\to C\in \RR,$ as $ h \to 0$.
\begin{lemma}\label{lemma:lt}
Let $ L_{h,N},\bar{L}_{h,N},\hat{L}_{h,N}$ be the local time estimators  defined above and in \eqref{def:LT:discr}. Then, 
\begin{align}
 \frac{ L_{h,N}}{hN} & \convas[N\to\infty]  \Linfty   \label{eq:cons:lt} \\
 \frac{ \bar{L}_{h,N}}{hN} &  \convas[N\to\infty]  \bar{L}(h) =\Linfty +\l_{1}\sqrt{h}+O(h) \label{eq:cons:lt-1} \\
 \frac{ \hat{L}_{h,N}}{hN} &  \convas[N\to\infty]   \hat{L}(h) = \Linfty
 +\l_{2}
 \sqrt{h}+
 O(h) \label{eq:cons:lt-2}
  \end{align}
  where $ \Linfty$ is in \eqref{eq:lt:limit}, $\bar{L}(h),\hat{L}(h)$ are functions of the parameters and of $h$,
  \[
\l_{1} =  \sqrt{\frac{
\pi}{2}}
 \dfrac{ |b_+| b_-}{\sigma_+ + \sigma_-}>0 
 \mbox{ and }
 \l_{2} = 
 -
 \frac{3\sqrt{\pi}}{4\sqrt{2}}
 \frac{|b_{+}|b_{-}}{|b_{+}|+b_{-}}
 \left( 
 \frac{|b_{+}|}{\sigma_{+}}
+
  \frac{b_{-}}{\sigma_{-}}
  +
\frac{|b_{+}|+b_{-} }{
\sigma_-+\sigma_+}
  \right)<0.
  \]
Moreover,
\begin{align*}
\sqrt{N}\left( \frac{ L_{h,N}}{hN} - \Linfty\right) \xrightarrow[N\to \infty]{\mathrm{law}}
 \cN(0, \Gamma_{\l}),
  \\
 \sqrt{N}\left(\frac{ \bar{L}_{h,N}}{hN} - \bar{L}(h) \right)\xrightarrow[N\to \infty]{\mathrm{law}}
 \cN(0, \Gamma_{\bar\l}),
\\
\sqrt{N}\left(  \frac{ \hat{L}_{h,N}}{hN} -  \hat{L}(h) \right)\xrightarrow[N\to \infty]{\mathrm{law}}
 \cN(0, \Gamma_{\hat\l}),
  \end{align*}
where $\Gamma_{\l},\Gamma_{\bar\l},\Gamma_{\hat\l}$ can be computed from \eqref{eq:def:gamma2}.
\end{lemma}
This lemma is proved in Section~\ref{sec:gme}.
The behavior of the three estimators is consistent with \eqref{eq:lt:limit} if we only look at the first term in the expansion of the limit. However, the $\sqrt{h}$ term is non zero for the limit of $\bar{L}_{h,N}$ and
 $\hat{L}_{h,N}$. So the only estimator having the same limit behavior  as the local time, when observations are in low frequency and the time horizon goes to infinity, is ${L}_{h,N}$. 
Furthermore, note that $L_{h,N}$ is the only one that does not require previous knowledge of the volatility parameters.

\subsection{Derivation of the estimators and comments on the results} \label{rm:derivation}

\paragraph{Derivation of the estimators.}
The dMLE in \eqref{eq:QMLE_dis_time} has been proposed in \citep{lp2} and is the maximum of a discretized classical likelihood for the solution of a SDE, based on the Girsanov transform. It can also be derived as a LSE with contrast function 
\[
	\sum_{k=0}^{N-1} \big|X_{(k+1)h}-X_{k h} - h b^+ \ind{ \{  X_{k h} \geq 0 \} } -  h b^{-} \ind{ \{  X_{k h} < 0 \} } \big|^{2}.
\]

\bigskip

The GME \eqref{eq:GME} is related to the one proposed for a different threshold SDE in~\citep{HuXi}. 
If the process is ergodic, the GME estimator is derived from the inversion of quantities such as $\EE\left[(X_{\infty})^{p}1_{\{\pm X_{\infty}\}>0}\right]$, where $p\in \NN$ and $X_{\infty}$ is a r.v.~whose distribution is the invariant one (it has density~\eqref{inv_dens}, see Section~\ref{app_fun_sys}).
In particular, one can compute a notion of asymptotic occupations time as
\begin{equation}\label{def:Q:stat}
{Q}_\infty^{\pm} 	
:=
\EE[\ind{\{ \pm X_\infty \geq 0 \}}]
=
\frac{|b_\mp|}{b_-+|b_+|},
\quad
{Q}_\infty^{\pm,1} 	
:=
\EE[ X_\infty\ind{\{ \pm X_\infty \geq 0 \}}]
=
\frac{\pm|b_\mp| \sigma_\pm^{2}}{2|b_\pm |(b_-+|b_+|)} . 
\end{equation}
These quantities can also be approximated from discrete observations using the ergodic theorem (cf.~next equations \eqref{erg:Q:LLN} and \eqref{erg:Q1:LLN}).
From the approximations $L_{h,N}$ to
$\Linfty $ (cf. next equations~\eqref{eq:lt:limit}-\eqref{eq:cons:lt}) and $\cQ^{\pm}_{h,N}$ to ${Q}_\infty^{\pm}$
 one can derive the estimators $ \overline{b}^\pm_{h,N}$ in \eqref{eq:GME}.
 Note that similar estimators could be defined using $\hat L_{h,N}$ or $\bar L_{h,N}$, but these would not have the correct limit (they would not be consistent) as $N\to\infty$, for fixed $h>0$, because of Lemma \ref{lemma:lt}. However, for small $h$ they could give good approximations.

Similarly, the definition of the volatility estimator in \eqref{eq:def:est:sigma} follows from the same approximations and the asymptotic value  $Q^{\pm,1}_{\infty}$ of $Q^{\pm,1}_{h,N}$,
corresponding to the first conditional moment used in \cite{HuXi}. Note the dependence on $\sigma_{\pm}$ in $Q^{\pm,1}_{\infty}$, not present in ${Q}_\infty^{\pm}$ and $\Linfty $.

The local time estimators are special cases of the class of statistics considered for instance in~\cite{j1} for Brownian diffusions for fixed time horizon.
In the context of threshold diffusions for fixed time horizon,
	${L}_{h,N}$ was considered in~\cite{lp2}, 
	the estimator related to the number of crossings, $\bar{L}_{h,N}$, was considerd in~\cite{port1,lmt1}, 
	and $\hat{L}_{h,N}$ was considered in~\cite{LP}.

\smallskip

\paragraph{On the equivalence of drift GME and dMLE, beyond the ergodic regime.}
In Lemma \ref{lemma:equivalence} we see that, in the ergodic case, GME and dMLE are equivalent in long time. This is because the dMLE is essentially given by the sum of two parts, one corresponding to the GME, the other given by the final and initial value of the process at $T$, normalized with the occupation time:
\begin{equation} \label{eq:diff:estim}
	\widehat{b}^\pm_{h,N}  = \overline{b}^\pm_{h,N} \pm \frac{\{X_{h N}\}^\pm -\{X_{0}\}^\pm}{Q^\pm_{h,N}} \qquad \PP\text{-a.s.}
\end{equation}
This equation is proved in the proof of Lemma~\ref{lemma:equivalence} in Section~\ref{sec:relation}. 
In the ergodic case, the fraction in~\eqref{eq:diff:estim} vanishes  
as $T\to \infty$ (see Lemma~\ref{lem:conv0} in Section~\ref{sec:relation}). 

A complete statistical analysis of the estimators from continuous time observations is proposed in~\cite[Propositions~3-7, 10-12]{lp2}, more precisely the estimators behavior is studied even if the process is not ergodic.
The results of the present document imply the same results for the dMLE in high frequency and infinite horizon. 
In this paper we deal only with the ergodic case in Theorem~\ref{th:joint:CLT}). 
The proof, for non ergodic cases, is the same as the one provided in Section~\ref{sec:mle} since Lemmas~\ref{lemma:control:diff:disc:cont}-\ref{lem:skew:density} also hold for non-ergodic cases. In the null recurrent case with non-vanishing drift, one needs to impose that  $h_N T_N^{3/2}\to 0$.

When the process is not ergodic, GME is not necessarily equivalent to the dMLE. 
For instance, in long time, in the transient case, the fraction in~\eqref{eq:diff:estim} does not vanish on the side where the process stays indefinitely, and it actually is the term providing the convergence of the estimators to the drift parameter. 
This is reminiscent of the drift MLE in a simple drifted Brownian motion (in finance, the Black-Scholes model), which does not depend on the intermediate observations, but only on the process value at $T$ and therefore, does not depend on the frequency of the observations. 
Note that if the drift is not $0$, the Black-Scholes model is transient, while the model in \eqref{eq:DOBM} can be transient or recurrent/ergodic depending on the drift parameters.
The normalized final value of the process at time $T$ is also the leading term for the model~\eqref{eq:DOBM}, when at least one of the drift parameters is such that the process is not mean reverting, which includes some null-recurrent cases.

\paragraph{Drift estimation in high frequency and infinite horizon.}

Condition~\eqref{cond:TN} is in line with what has been obtained in the literature so far for non-threshold diffusions, and it is the best condition obtained so far for threshold ones.
Up to the authors knowledge, Theorem~\ref{th:joint:CLT} is the first result for threshold diffusions allowing for deterministic initial conditions (instead of starting at the stationary distrubutions) when looking at high-frequency observations over an infinite time horizon. 

Theorem~\ref{th:joint:CLT} holds also for non-equally spaced high frequency observations. More precisely, the process is observed at times $0=t_0<t_1<\ldots<t_N=T_N$ and $h_N$ represents the maximal observation lag: $h_N=\sup_{k=0,\ldots, N-1} (t_{k+1}-t_{k})$.


\paragraph{Multi-threshold.}

The method we are presenting here should also work in a multi-threshold setting, at least in the high frequency and long time setting, where it works analogously to \cite[Supplementary material Section S.3]{mazzonetto2020drift}. 
In the GME case, with more than one threshold, the definition of the estimators would involve more complicated, implicit expressions, that one would need to solve using numerical methods, while, if there is only one threshold, explicit expressions are available. Therefore, for the sake of simplicity of the exposition and ease of implementation, we stick here to the case with only one threshold.


\section{Numerical studies} \label{sec:numerics}

\subsection{Comparison of drift estimators}
\begin{figure}[ht!]
\centering
	\includegraphics[width=\textwidth]{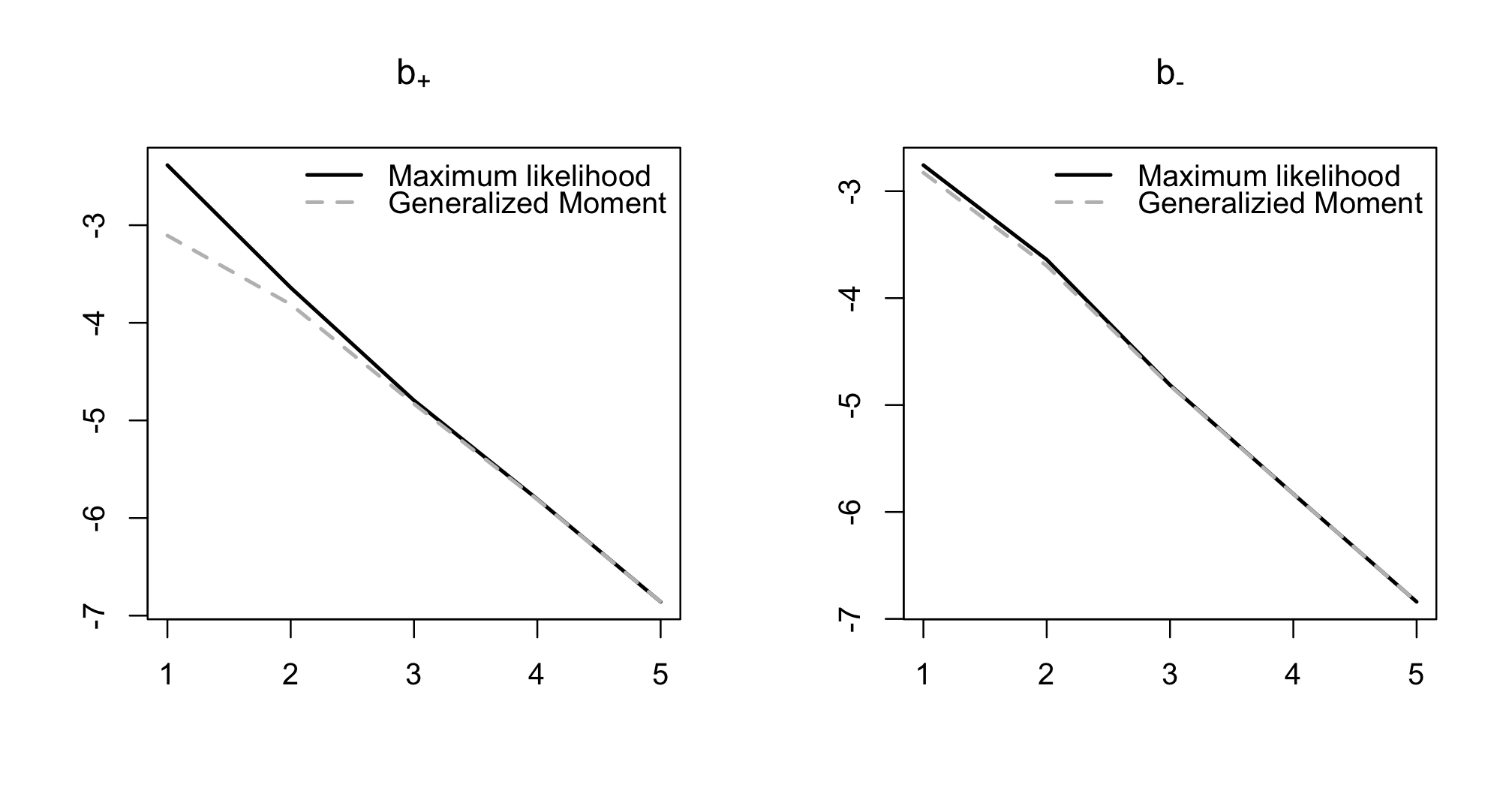}
\caption{We compare the GME $ \overline{b}^\pm_{h,N}$ and the dMLE  $\hat{b}^\pm_{h,N}$, with $h=1$ fixed, varying the time horizon (length of the time series) $N$. We plot in this figure the mean squared error of each estimator, on $10^{3}$ simulated time series, as a function of the length of the time series, in a log-log plot in base $10$.
}
\label{fig:mse:drift}
\end{figure}

We compare here estimators $ \overline{b}^\pm_{h,N}$ with  $\hat{b}^\pm_{h,N}$, with $h>0$ fixed, as the number of observations goes to infinity, plotting the mean squared error of the estimators (MSE), for a process with parameters as in Table \ref{table:simulation_parameters}, simulated via the Euler-Maruyama scheme. 
\begin{table}[htp]
\centering
\begin{tabular}
{c|c|c|c}%
 $b_-$ & $b_+$ & 
$\sigma_-$ & $\sigma_+$ 
  \\ \hline
   $0.02$ & $-0.01$ & 
$0.07$ & $0.10$ 
\end{tabular}
\caption{Simulations parameters.
\label{table:simulation_parameters}}
\end{table}
Clearly, because of Theorems \ref{th:fixed:h} and \ref{th:joint:CLT}, the asymptotic behavior of the estimators is the same. However, our numerical experiments suggest that the GME slightly outperforms the dMLE, not by a significant margin, but consistently over different sets of parameters and time horizons (see Figure \ref{fig:mse:drift}). This can be explained by the discussion in Section \ref{rm:derivation}, since we are in the ergodic setting.

\subsection{The asymptotic normality constant in the drift estimator}\label{sec:num:drift}

\begin{figure}[ht!]
\centering
	\includegraphics[width=\textwidth]{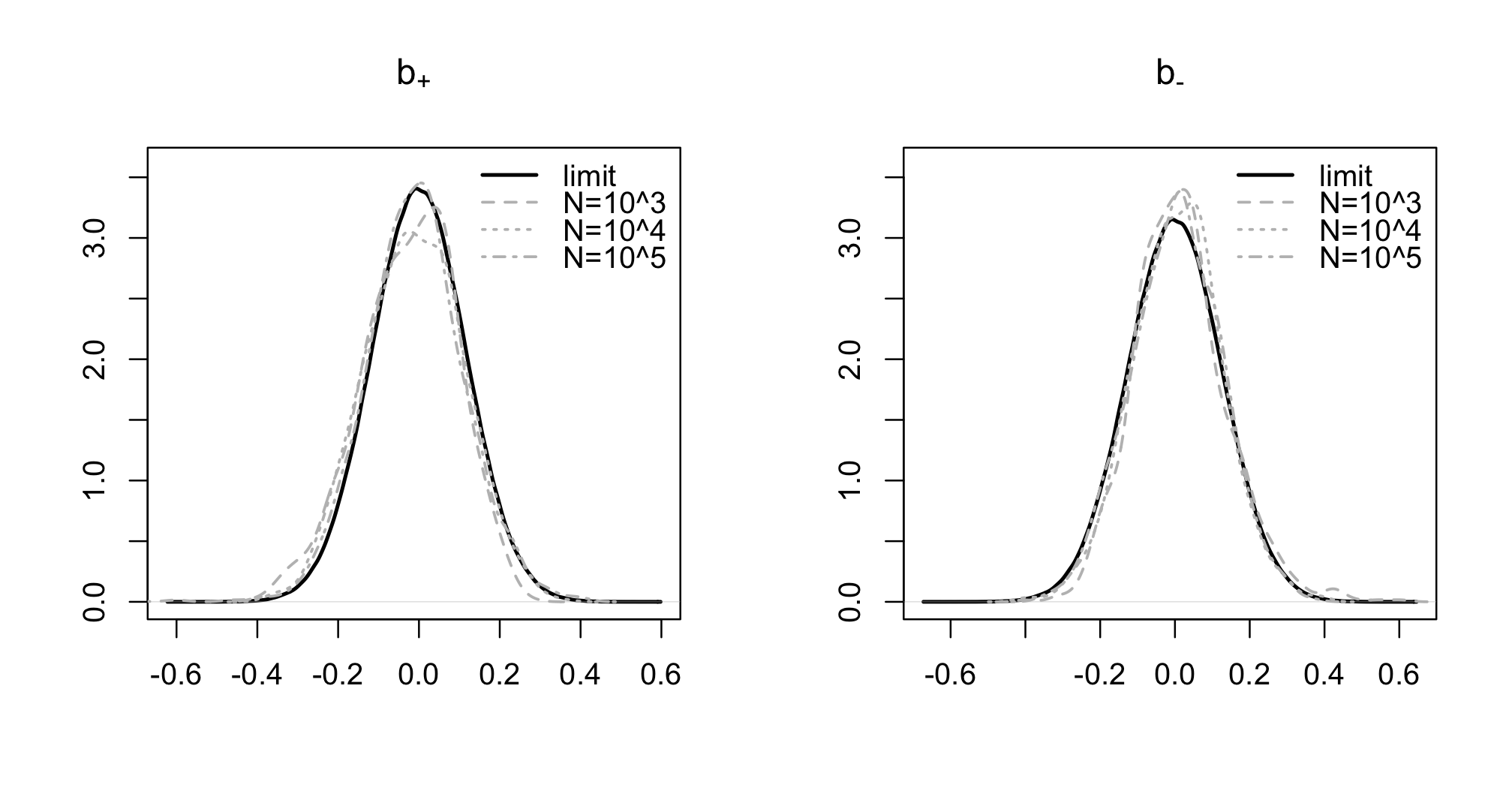}
\caption{
CLT in Theorem~\ref{th:fixed:h}-ii), with parameters as in Table~\ref{table:simulation_parameters}. We plot the density of the theoretical distribution of the estimation error, with variance estimated on $10^{3}$ simulated paths, and compare it with
the distribution of the rescaled error on $n=10^3$ paths, and $N=10^3,10^4,10^5$ observations on each path.
\label{fig:clt:drift}}
\end{figure}

We look now at the CLT in Theorem \ref{th:fixed:h}
- ii). The limit covariance matrix is given in \eqref{eq:def:gammaq}, for which we do not have an analytical expression. However, it can be evaluated numerically, on synthetic data, simulating several realizations of the process and then computing the empirical covariance. We do so, and plot in Figure \ref{fig:clt:drift} the empirical distribution of the error, rescaled with $\sqrt{N}$, for various lengths of the observed time series, generated with parameters as in Table \ref{table:simulation_parameters} via the Euler-Maruyama scheme. A similar numerical representation for Theorem \ref{th:joint:CLT} can be found in \cite{lp2}, where the continuous time equivalent of Theorem \ref{th:joint:CLT} is stated and checked numerically.

In view of the independence of the limit Gaussians in Theorem \ref{th:joint:CLT}, an interesting question is whether the limit cross covariance vanishes or not also in Theorem \ref{th:fixed:h}. In our simulations, we can never reject the null hypothesis of $0$ limit cross covariance. However, this numerical evidence is not conclusive since the limit cross covariance could be small but non-zero, and still be consistent with our simulations.

These simulations provide a numerical validation of  Theorem \ref{th:fixed:h}. However, this approach generally cannot be used to estimate the error when dealing with parameters inference from empirical data, since in this case estimation can only be based on one observed time series and several realizations are not available. Estimating the limit covariance can be necessary, for example in order to perform a test on the estimated parameters. In this case, one way to estimate the limit covariance is to split the time series in several batches and use them as separate realizations, as explained in \cite[Section 1.10.1]{brooks2011handbook}. Alternative methods are presented e.g. in \cite[Section 1.10.2 and 1.10.3]{brooks2011handbook}, see also 
\cite[Section 1.12.9]{brooks2011handbook} for subsampling techniques.


\section{Technical material and proofs} \label{sec:proofs}

We collect in this section all the mathematical proofs. 
In Section~\ref{app_fun_sys} we give some preliminary technical results related to the law of the process via its infinitesimal generator. 
In Section~\ref{sec:relation} we prove the asymptotic equivalence of the GME and dMLE and Theorem~\ref{th:disc}.
In Section \ref{sec:gme}, we provide the proofs related to the GME, with fixed time lag and long time horizon. 
Finally, in Section \ref{sec:mle} the proofs related to the dMLE, in high frequency and with long time horizon. 

\subsection{Technical tools: scale function, speed measure, fundamental system and resolvent kernel}\label{app_fun_sys}

The infinitesimal generator $\mathcal{A}$ of the process $X$ solution to \eqref{eq:DOBM} can be written as
\begin{gather*}
    \mathcal{A} f=\frac{1}{2}\sigma^2(x)e^{-h(x)}\frac{\dd}{\dd x}\left(
    e^{h(x)}\frac{\dd f(x)}{\dd x}\right)\ 
    \text{ with }
    h(x)=\int_0^x \frac{2b(y)}{\sigma^2(y)}\vd y,
\end{gather*}
cf.~e.g.~\citep{ito} for a general discussion and \citep{lp2} for the diffusion in \eqref{eq:DOBM}. 
The process $X$ is a diffusion 
(roughly-speaking, a strong Markov process with continuous sample paths)
which is fully characterized by its \emph{speed measure} $M$
with density $m$ and \emph{scale function} $S$ given by
\begin{equation}
    \label{eq:speend-n-scale}
    m(x)\eqdef\frac{2}{\sigma(x)^2}\exp( h(x) )
\text{ and }
S(x)\eqdef\int_0^x \exp( -h(y) )\vd y.
\end{equation}
When the process is positive or null recurrent, the speed measure is a stationary measure (invariant measure for the transition semigroup).
When $b_+<0$, $b_->0$ we say that the process is \emph{ergodic} and its \emph{invariant distribution} is the re-normalized speed measure, denoted by $\mu$, with density given by
\begin{equation}\label{inv_dens}
    \frac{m(x)}{M(\RR)}
    =\begin{cases}
	\dfrac{ 2}{\sigma_+^2}\times\dfrac{|b_+| b_-}{b_-+|b_+|}e^{\frac{-2|b_+| x }{\sigma_+^2}}&\text{ if }x\geq0,\\
		\dfrac{ 2}{\sigma_-^2}\times
		\dfrac{b_-|b_+|}{b_-+|b_+|}e^{\frac{2b_- x}{\sigma_-^2}}&\text{ if }x<0.
    \end{cases}
\end{equation}
Here $M(\RR)=\int_\RR m(x)dx = (|b_+|+b_-)/(|b_+|b_-)$. {(Note that there is a typo in \citep{lp2}.)}
Throughout the paper, we write $X_\infty$ for a r.v. independent of $X,$ following the stationary distribution, i.e., with density \eqref{inv_dens}. 
Let us consider equation
\begin{equation}\label{eqn:heat}
(\mathcal{A} -\lambda)u= 0.
\end{equation}
For fixed $\lambda>0$,  the set of solutions to \eqref{eqn:heat}
is a two-dimensional vector space. There exist two continuous, positive functions $\psi_\lambda, \varphi_\lambda$ solution to \eqref{eqn:heat}, with $\psi_\lambda$ increasing from $0$ to $\infty$ and $\varphi_\lambda$ decreasing from $\infty$ to $0$, such that $\psi_\lambda(0)=\varphi_\lambda(0)=1$, called the minimal functions. Such functions can be taken as a basis for the space of solutions to \eqref{eqn:heat}.
With coefficients as in \eqref{sigmaDOBM},  the minimal functions are explicit and given as
\begin{align}
    \label{increasingsolution}
    \psi_\lambda(x)&=\begin{cases}
    \exp\left(x\frac{-b_-+\sqrt{b_-^2+2\sigma_-^2\lambda}}{\sigma_-^2}\right)&\text{ if }x<0\\
    \kappa_+ \exp\left(x\frac{-b_++\sqrt{b_+^2+2\sigma_+^2\lambda}}{\sigma_+^2}\right)
    + \delta_+ \exp\left(x\frac{-b_+-\sqrt{b_+^2+2\sigma_+^2\lambda}}{\sigma_+^2}\right)&\text{ if }x\geq0,
\end{cases}
    \\
    \label{decreasingsolution}
    \varphi_\lambda(x)&=
    \begin{cases}
	\kappa_-\exp\left(x\frac{-b_--\sqrt{b_-^2+2\sigma_-^2\lambda}}{\sigma_-^2}\right)
	+\delta_-\exp\left(x\frac{-b_-+\sqrt{b_-^2+2\sigma_-^2\lambda}}{\sigma_-^2}\right)&\text{ if }x<0,\\
	\exp\left(x\frac{-b_+-\sqrt{b_+^2+2\sigma_+^2\lambda}}{\sigma_+^2}\right)&\text{ if }x\geq0
    \end{cases}
\end{align}
with 
\begin{align*}
    \kappa_+&\eqdef
\frac{-b_-\sigma_+^2 + b_+\sigma_-^2 + \sigma_-^2\sqrt{b_+^2 + 2\lambda\sigma_+^2} + \sigma_+^2\sqrt{b_-^2 + 2\lambda\sigma_-^2}}{2\sigma_-^2\sqrt{ b_+^2 + 2\lambda\sigma_+^2 }},\\
\delta_+&\eqdef
\frac{b_-\sigma_+^2 - b_+\sigma_-^2 + \sigma_-^2\sqrt{b_+^2 + 2\lambda\sigma_+^2} - \sigma_+^2\sqrt{b_-^2 + 2\lambda\sigma_-^2}}{2\sigma_-^2\sqrt{b_+^2 + 2\lambda\sigma_+^2}},\\
\kappa_-&\eqdef
\frac{-b_-\sigma_+^2 + b_+\sigma_-^2 - \sigma_-^2\sqrt{b_+^2 + 2\lambda\sigma_+^2} + \sigma_+^2\sqrt{b_-^2 + 2\lambda\sigma_-^2}}{2\sigma_+^2\sqrt{b_-^2 + 2\lambda\sigma_-^2}},\\
\delta_-&\eqdef
\frac{b_-\sigma_+^2 - b_+\sigma_-^2 + \sigma_-^2\sqrt{b_+^2 + 2\lambda\sigma_+^2} + \sigma_+^2\sqrt{ b_-^2 + 2\lambda\sigma_-^2 }}{2\sigma_+^2\sqrt{b_-^2 + 2\lambda\sigma_-^2}}.
\end{align*}

Using the minimal functions, one can define the
 Wronskian of the diffusion as
\[
\begin{split}
W_\lambda
& =
\varphi_\lambda(x) \frac{\partial_x \psi_\lambda(x)}{ \partial_x S(x)}
-
\psi_\lambda(x) \frac{\partial_x \varphi_\lambda(x) }{ \partial_x S(x)} \\
&=\frac{-b_-+\sqrt{b_-^2+2\sigma_-^2\lambda}}{\sigma_-^2}
-
\frac{-b_+-\sqrt{b_+^2+2\sigma_+^2\lambda}}{\sigma_+^2}\\
&=\frac{-b_-+\sqrt{b_-^2+2\sigma_-^2\lambda}}{\sigma_-^2}+
\frac{b_+ +\sqrt{b_+^2+2\sigma_+^2\lambda}}{\sigma_+^2}
\end{split}
\]
where $S$ is the scale function in \eqref{eq:speend-n-scale}. Note that
\begin{equation}\label{eq:as:W}
W_\lambda = \lambda\frac{b_-+|b_+|}{b_-|b_+|} +O(\lambda^2)
\end{equation}
Let us define the \emph{resolvent kernel} 
of \eqref{eq:DOBM} as the Laplace transform of the transition density
\begin{equation}\label{resolvent}
r(\lambda,x,y)=\mathcal{L}_\lambda (p(\cdot,x,y))=\int_0^\infty p(t,x,y)\exp(-\lambda t) dt.
\end{equation}
It can be computed using the speed measure \eqref{eq:speend-n-scale}, the minimal functions \eqref{increasingsolution}, \eqref{decreasingsolution} and the Wronskian, as
\begin{equation}\label{resolventmin}
r(\lambda,x,y)=\frac{m(y)}{W_\lambda}
\begin{cases}
 \psi_\lambda(x)\varphi_\lambda(y) \quad \mbox{ if }x<y, \\
 \varphi_\lambda(x) \psi_\lambda(y) \quad \mbox{ if }x>y.
\end{cases}
\end{equation}
We can use it as follows to compute Laplace transforms of expectations of functions $f(\cdot)$ of $X_t$ starting from a point $z$, if we can exchange the order of time integral and expectation, since
\begin{equation}\label{eq:propr:res}
R_\lambda f(z)
:=
\cL_\lambda \EE_z [f(z,X_\cdot)]
=
\int_0^\infty  \EE_z [e^{-\lambda t} f(z,X_t) dt ]
=
\int_\RR r(\lambda,z,y)f(z,y)dy.
\end{equation}

\subsection{The relation between different estimators}\label{sec:relation}

\begin{lemma}
\label{lem:conv0}
Let $X$ be solution to \eqref{eq:DOBM}. Then
\begin{equation}
t^{-p} X_t \convp[t \to \infty] 0, \mbox{ for any } p \in(0,+\infty).
\end{equation}
\end{lemma}
\begin{proof}

\textit{Step 1}: We first prove that the statement holds with initial condition $X_0=0$. 
We show that
\begin{equation}
 \EE[|X_t|] \to c\in \RR, \mbox { for } t\to\infty 
\end{equation}
which implies the statement. 
We use that
$\EE[|X_t|]=\EE[X_t^+]+\EE[X_t^-]$ and write
\[
\cL_\lambda \EE[X_\cdot^+] =  \int_0^\infty r(\lambda, 0, y) y dy
\]
with $r$ in \eqref{resolvent}.
Using  \eqref{resolventmin}, we compute
\[
 	\cL_\lambda \EE[X_\cdot^+]  
	=  \frac{\psi_\lambda(0)}{W_\lambda}\int_0^\infty m(y)\varphi_\lambda(y) y dy
	 =
 	\frac{1}{W_\lambda}\int_0^\infty m(y) \varphi_\lambda(y) y dy.
\]
We have
\[
\begin{split}
\int_0^\infty m(y) \varphi_\lambda(y) y dy
&=
\frac{2}{\sigma_+^2} \int_0^\infty 
\exp\left(y\frac{b_+-\sqrt{b_+^2+2\sigma_+^2\lambda}}{\sigma_+^2}\right) y dy\\
&= 
\frac{2\sigma_+^2}{(b_+-\sqrt{b_+^2+2\sigma_+^2\lambda})^2}
=
\frac{
(\sqrt{b_+^2+2\sigma_+^2\lambda}+b_+)^2
}{2\sigma_+^2 \lambda^2}
\end{split}
\]
and, usign also \eqref{eq:as:W}, as $\lambda \to 0$, 
\[
\cL_\lambda \EE[X_\cdot^+] = \frac{\sigma_+^2 b_-}{2|b_+| (b_-+|b_+|)}\frac{1}{\lambda }+O(1).
\]
So, using the final value theorem (Tauberian theorem), we see that as $t \to \infty$
\[
\EE[X_t^+] \to \frac{\sigma_+^2 b_-}{2|b_+| (b_-+|b_+|)}
\]
which is $\EE[X_\infty^+]$ from \eqref{inv_dens}.
With a similar computation we prove an analogous result for $\EE[X_t^-]$, and therefore as $t \to\infty$ we obtain $\EE[X_t] \to  \EE[X_\infty]$ and $\EE[|X_t|] \to  \EE[|X_\infty|]$, where the limits can be directly computed and are finite. Markov's inequality completes the proof of this first step.

\textit{Step 2}: The statement for $X_0=0$ implies the statement for any $X_0\in \RR$, say $X_0=x$. Let us denote the process with such initial condition by $X$.
Let $\tau$ be the stopping time identifying the first time for which $X_t=0$. Note that $\tau$ is a.s.~finite since $X$ is recurrent.
The process $Y_{t}:=X_{t+\tau}$ is a weak solution to~\eqref{eq:DOBM} with $Y_0=0$. Hence, by the strong Markov property, we deduce that $\EE[|Y_t|]$ converges to a constant as $t$ goes to infinity.
To prove the statement it suffices to show that  $\sup_{t>0} \EE[|Y_{t}-X_t|]$ is bounded by a constant.
From \eqref{eq:DOBM}, follows that $X_{t+\tau}-X_t =\int_t^{t+\tau} b(X_s) \vd s + \int_t^{t+\tau} \sigma(X_s) \vd W_s$.
Then, H\"older's inequality and It\^o isometry yield for all $t>0$
\[ \EE[|Y_{t}-X_t|]  = \EE[|X_{t+\tau}-X_t|] \leq 2 \max\{|b_+|, b_-,\sigma_-,\sigma_+\} (1+\EE[\tau]),\]
and $\EE[\tau]<\infty$ because the process is positive recurrent. The proof is thus completed.
\end{proof}

\begin{lemma}\label{lemma:equivalence}
For $h>0$ fixed 
\[
	\sqrt{N} \left(\overline{b}^\pm_{h,N}- \widehat{b}^\pm_{h,N} \right) \convp[N\to \infty] 0.
\]
For $h_N=T_N/N$, with $T_N$ and $N$ satisfying condition~\eqref{cond:TN} and such that $\lim_{N\to\infty} T_N =\infty$,
\[
	\sqrt{h_N N} \left(\overline{b}^\pm_{h_N,N}- \widehat{b}^\pm_{h_N,N} \right) \convp[N\to \infty] 0.
\]
For $T>0$ fixed, $h_N=T/N$
as $N\to \infty$,
\[
\overline{b}^\pm_{h_N,N}- \widehat{b}^\pm_{h_N,N}\convp[N\to\infty]
	\mp \frac{\{X_{T}\}^\pm -\{X_{0}\}^\pm}{Q^\pm_{T}}.
\]
\end{lemma}
\begin{proof}
It\^o-Tanaka formula
establishes that the following equalities hold $\PP$-a.s.:
\begin{equation} \label{Tanaka:cont}
\begin{split}
	\{X_T\}^\pm -  \{X_0\}^\pm 
	& =
		\pm \cM^{\pm}_T
+ \frac{1}{2} L_T(X),
		\end{split}
\end{equation}
Moreover, it is not difficult to see, similarly to~\cite[(S2.6) in supplementary material]{mazzonetto2020drift}, that for any $h,N$, we have $\PP$-a.s. 
\begin{equation} \label{Tanaka:disc}
\begin{split}
	\{X_{hN}\}^\pm -  \{X_0\}^\pm 
	& =
		\pm \cM^{\pm}_{h,N}
+ \frac{1}{2} L_{h,N}.
\end{split}
\end{equation}
Hence, $\PP$-a.s.
\[
	\sqrt{N} \left(\overline{b}^\pm_{h,N}- \widehat{b}^\pm_{h,N} \right)= 
	\mp \frac{1}{\sqrt{h}}  \frac{h N}{\sqrt{h N}}  \left( \frac{\{X_{h N}\}^\pm -\{X_{0}\}^\pm}{Q^\pm_{h,N}}\right).
\]
For $h$ fixed, when ${N\to\infty}$, we have $ Q^\pm_{h,N}/(h N)\convas Q^\pm_{\infty}$, because of \eqref{erg:Q:LLN}.
For $h_N=T_N/N$ and $T_N=h_N N$ such that $\lim_{N\to \infty} T_N =\infty$ and \eqref{cond:TN}, we have $ Q^\pm_{h_N,N}/T_N\convas {Q}^\pm_{\infty}$ when ${N\to\infty}$, because of the ergodic theorem and
Lemma~\ref{lemma:control:diff:disc:cont}.
Lemma~\ref{lem:conv0} implies that $ X_t/\sqrt{t} \convp 0$ when ${t\to\infty}$. The first two statements follow.

Now, when $T>0$ fixed, $h_N=T/N$, we have $
Q^\pm_{h_N,N}\convas  Q^\pm_T$ when $N\to\infty$ and the last statement follows.
\end{proof}

\subsubsection{Proof of Theorem~\ref{th:disc}}\label{pr:thdisc}
The result states that, for fixed time horizon, in high frequency, the drift dMLE from discrete observations converges in probability towards the drift MLE from continuous observations~\eqref{eq:QMLE_ct} and provides the convergence rates.
It also states that 
the generalized moment estimator \eqref{eq:GME}, converges for fixed time and in high frequency to its continuous-time analog
\[
\overline \ell_T := \begin{pmatrix}
	-	L_T 
/(2Q^{+}_T) 
, & {	L_T / 
(2Q^{-}_T )} 
	\end{pmatrix}.\]
In Lemma~\ref{lemma:equivalence} it is shown that, for fixed time and in high frequency, 
\[
\widehat{b}^\pm_{h_N,N} - \overline{b}^\pm_{h_N,N} \convp[N\to\infty]
	\pm \frac{\{X_{T}\}^\pm -\{X_{0}\}^\pm}{\cQ^\pm_{T}}
\]
and, from~\eqref{Tanaka:cont}-\eqref{Tanaka:disc} in the proof of Lemma~\ref{lemma:equivalence}, it holds $\PP$-a.s.~that
\[ 
	\widehat b^\pm_{h_N,N} - \beta^\pm_T + \overline \ell^\pm_T - \overline b^\pm_{h_N,N} 
	= \pm (\{X_T\}^\pm -\{X_0\}^\pm)\frac{\cQ^\pm_{T}-\cQ^\pm_{h_N,N}}{\cQ_T^\pm \cQ^\pm_{h_N,N}}.
\]
The speed of convergence of the occupation time is at least $\sqrt{N}$ (see~\cite{LP}) which is larger than $N^{1/4}$.
Hence, $\PP$-a.s.~
\[ 
 \overline b^\pm_{h_N,N}
-\overline \ell^\pm_T =
\widehat b^\pm_{h_N,N} - \beta^\pm_T +O(1/\sqrt{N}).
\]
We are thus reduced to prove the result for $\widehat{b}^+_{h_N,N}$. 
Then, \eqref{Tanaka:cont} and \eqref{Tanaka:disc} ensure that $\PP$-a.s.:
\[
	\widehat b^\pm_{h_N,N} - \beta^\pm_T =\mp\frac{L_{h_N,N} - L_T}{2 \cQ^\pm_T} + \cM^\pm_T \frac{\cQ^\pm_T - \cQ^\pm_{h_N,T}}{\cQ^\pm_T \cQ^\pm_{h_N,T}}
\]
Recall $h_N=T/N$.
The consistency, i.e.~$(\widehat{b}^+_{h_N,N} -  \beta^+_{T}, \widehat{b}^-_{h_N,N} - \beta^-_{T} ) \xrightarrow{\PP} 0$ when $N\to \infty$, was already provided in~\cite[Lemma~1]{lp2}, and the convergence rate is based on the one of the local time approximation~\cite[Proposition~2]{mazzonetto2019rates}. 
The proof is thus completed.

\subsection{Proofs of asymptotic behaviors in low frequency and long time: GME in Theorem~\ref{th:fixed:h} and Theorem~\ref{th:sigma:h} and local time in Lemma~\ref{lemma:lt}}
\label{sec:gme}

Let us define the 2-dimensional discrete time process $(\X_k^h)_{k\in\NN}=(X_{hk},X_{h(k+1)})_{k\in\NN}$. In the ergodic case, we set $\X_\infty^h=(X_\infty,X_\infty^h)$, a r.v.~that follows its stationary distribution, where $X_{\infty}$ is a r.v.~following the stationary distribution $\mu$ of \eqref{inv_dens} and $X_{\infty}^h$ is its evolution over a time $h$ following equation \eqref{eq:DOBM}. 

We also define $\tilde{X}= (\tilde{X}_t)_{t}$ a solution to \eqref{eq:DOBM} with $\tilde{X}_{0}\eqlaw X_\infty$ and $(\tilde{\X}_k^h)_{k\in\NN}=(\tilde{X}_{hk},\tilde{X}_{h(k+1)})_{k\in\NN}$.

\begin{lemma}\label{lemma:ergodic}
Let $h>0$ fixed, $\pm b_\pm <0$. 
\begin{enumerate}[I]
\item
The discrete time process $(X_{hk})_{k\in\NN}$ is ergodic, with invariant measure in \eqref{inv_dens}.  For any measurable function $f:\RR\to\RR^{n}$ s.t. $\EE[|f(X_{\infty})|]<\infty$,
\[
\frac{1}{N} \sum_{k=1}^{N}f(X_{hk})\convas[N\to\infty]\EE[f(X_{\infty})].
\]
\item
Moreover, if $|f(x)|\leq C |x|^{p}$ for some $C_{0},C_{1},p>0$ we have that
\[
\sqrt{N}\left(\frac{1}{N} \sum_{k=1}^{N}f(X_{hk})-\EE[f(X_{\infty})]\right)\convl[N\to\infty]  \cN(0,\gamma)
 \]
 where setting $\bar{f}(\cdot)=f(\cdot)-\EE[f(X_{\infty})]$ we have
\begin{equation}\label{eq:def:gamma}
\begin{split}
\gamma=\EE[\bar{f}(\tilde{X}_0)\bar{f}^{T}(\tilde{X}_0)]+2\sum_{k=1}^\infty  \EE[\bar{f}(\tilde{X}_0)\bar{f}^{T}(\tilde{X}_k)].
\end{split}
\end{equation}
\item
The process $(\X_k^h)_{k\in\NN}$ is ergodic and, for any measurable function $f:\RR^{2}\to\RR^{n}$ s.t. $\EE[|f(\X^{h}_{\infty})|]<\infty$, 
\[
\frac{1}{N} \sum_{k=1}^{N}f(\X^{h}_{k}) \convas[N\to\infty] \EE[f(\X^{h}_{\infty})].
\]
\item
Moreover, if $|f(x)|\leq C_{0}+C_{1} |x|^{p}$ for some $C_{0},C_{1},p>0$,
we have that
\[
 \sqrt{N}\left(\frac{1}{N} \sum_{k=1}^{N}f(\X_{h}^{k})-\EE[f(\X^{h}_{\infty})]\right) \convl[N\to\infty] \cN (0,\Gamma)
 \]
 where setting $\bar{f}(\cdot)=f(\cdot)-\EE[f(X_{\infty})]$ we have
\begin{equation}\label{eq:def:gamma2}
\begin{split}
\Gamma=\EE[\bar{f}(\tilde{\X}^{h}_0) \bar{f}^{T}(\tilde{\X}^{h}_0) ]+\sum_{k=1}^\infty  \EE[\bar{f}(\tilde{\X}^{h}_0)\bar{f}^{T}(\tilde{\X}^{h}_k)]
+
\EE[\bar{f}(\tilde{\X}^{h}_k)\bar{f}^{T}(\tilde{\X}^{h}_0)].
\end{split}
\end{equation}
\end{enumerate}
\end{lemma}

\begin{remark}
The quantities defined in \eqref{eq:def:gamma} and \eqref{eq:def:gamma2} are always non-negative, so \eqref{eq:def:gamma} and \eqref{eq:def:gamma2} are good definitions of variances in $[0,+\infty)$.
Moreover, these are real CLTs with convergence rate $\sqrt{N}$ if  $\gamma,\Gamma$ are postive definite.
\end{remark}
\begin{proof} The proof of I and III is completely analogous to \cite[Theorem 2]{HuXi}, the main tool being an application of the ergodic theorem in
\cite[Theorem
17.0.1]{MeynTweedie}. Let us now discuss IV. Again we start from 
\cite[Theorem
17.0.1]{MeynTweedie}, which only applies to uni-dimensional $f$.
The fact that this holds for a multidimensional $f$ follows using the Cramer-Wald theorem as described in \cite[Section 1.8]{brooks2011handbook}. Clearly it is enough to prove it for $n=2$. Let $(U_{1},U_{2})$ be a Gaussian r.v. in $\RR^{2}$ with covariance
{\small
\[
\EE\begin{pmatrix}
\bar{f}_{1}^2(\tilde{\X}^{h}_0)]+2\sum_{k=1}^\infty   [\bar{f}_{1}(\tilde{\X}^{h}_0)\bar{f}_{1}(\tilde{\X}^{h}_k)
&
 \bar{f}_{1}\bar{f}_{2}(\tilde{\X}^{h}_0)+
\sum_{k=1}^\infty   \bar{f}_{1}(\tilde{\X}^{h}_0)\bar{f}_{2}(\tilde{\X}^{h}_k)+\bar{f}_{2}(\tilde{\X}^{h}_0)\bar{f}_{1}(\tilde{\X}^{h}_k)
\\
 \bar{f}_{1}\bar{f}_{2}(\tilde{\X}^{h}_0)+
\sum_{k=1}^\infty   \bar{f}_{1}(\tilde{\X}^{h}_0)\bar{f}_{2}(\tilde{\X}^{h}_k)+\bar{f}_{2}(\tilde{\X}^{h}_0)\bar{f}_{1}(\tilde{\X}^{h}_k)
&
 \bar{f}_{2}^2(\tilde{\X}^{h}_0)+2\sum_{k=1}^\infty   \bar{f}_{2}(\tilde{\X}^{h}_0)\bar{f}_{2}(\tilde{\X}^{h}_k)
\end{pmatrix}.
\]}
For $\lambda_{1},\lambda_{2}\in \RR$ we have, from the univariate version as in \cite[Theorem~17.0.1]{MeynTweedie}, after some computations, that in law
\[
\lim_{N\to \infty}\frac{1}{ \sqrt{N}} \sum_{k=1}^{N}
(\lambda_{1}\bar{f}_{1}(\X_{h}^{k}) + \lambda_{2}\bar{f}_{2}(\X_{h}^{k}) )
= \lambda_{1}U_{1}+\lambda_{2} U_{2},
 \]
 because $\lambda_{1}U_{1}+\lambda_{2} U_{2}$ is a centered Gaussian with variance given by
 \[
 \begin{split}
&\lambda_{1} \left(\EE[\bar{f}_{1}^2(\tilde{\X}^{h}_0)]+2\sum_{k=1}^\infty  \EE[\bar{f}_{1}(\tilde{\X}^{h}_0)\bar{f}_{1}(\tilde{\X}^{h}_k)]\right)
+\lambda_{2}\left(\EE[\bar{f}_{2}^2(\tilde{\X}^{h}_0)]+2\sum_{k=1}^\infty  \EE[\bar{f}_{2}(\tilde{\X}^{h}_0)\bar{f}_{2}(\tilde{\X}^{h}_k)]\right)
\\
&+2\lambda_{1}\lambda_{2}
\left(
\EE[\bar{f}_{1}\bar{f}_{2}(\tilde{\X}^{h}_0)]+
\sum_{k=1}^\infty  \left(\EE[\bar{f}_{1}(\tilde{\X}^{h}_0)\bar{f}_{2}(\tilde{\X}^{h}_k)] 
+
\EE[\bar{f}_{2}(\tilde{\X}^{h}_0)\bar{f}_{1}(\tilde{\X}^{h}_k)]\right)
\right).
\end{split}
\]
Therefore Cramer Wald theorem implies \eqref{eq:def:gamma2} in the multivariate sense. The proof of II is analogous, with the addition of the fact that 
\[
\EE[\bar{f}_{1}(\tilde{X}^{h}_0)\bar{f}_{2}(\tilde{X}^{h}_k)] 
=
\EE[\bar{f}_{2}(\tilde{X}^{h}_0)\bar{f}_{1}(\tilde{X}^{h}_k)].
\]
The proof is thus completed.
\end{proof}
We prove that the convergence results in  Theorem \ref{th:fixed:h} hold for $\overline{b}^\pm_{h,N}$. The fact that they hold for $\widehat{b}^\pm_{h,N}$ as well follows from Lemma \ref{lemma:equivalence}. Let us recall the expressions of ${Q}_\infty^{\pm}$ and ${Q}_\infty^{\pm,1}$ in~\eqref{def:Q:stat}.

\begin{lemma}
For fixed $h>0$, let $Q^{\pm}_{h,N}$ as in \eqref{def:M:Q:disc}.
Then
\begin{equation}\label{erg:Q:LLN}
\frac{
 Q^{\pm}_{h,N}}{N}
\convas[N\to\infty] h {Q}_\infty^{\pm}=\frac{h |b_\mp|}{|b_+|+b_-} 
\end{equation}
and
\begin{equation}\label{erg:Q:CLT}
\sqrt{N}\left(
\frac{
 Q^{\pm}_{h,N}}{N} - \frac{h |b_\mp|}{|b_+|+b_-}\right) 
 \convl[N\to\infty]
 \cN(0,h^{2}\xi_\pm),
\end{equation}
with
\begin{equation}\label{def:gammapm}
\begin{split}
\xi_\pm:= \EE[(\ind{\{ \pm X_0 \geq 0\}}-{Q}_\infty^\pm)^2]+ 2\sum_{k=1}^\infty \EE[(\ind{\{ \pm X_0 \geq 0\}}-
{Q}_\infty^\pm)( \ind{\{ \pm X_k \geq 0\}}-{Q}_\infty^\pm) ].
\end{split}
\end{equation}
Let now $Q^{\pm,1}_{h,N}$ be the estimator in \eqref{def:Q1:disc}. Then
\begin{equation}\label{erg:Q1:LLN}
\frac{
 Q^{\pm,1}_{h,N}}{N}
 \convas[N\to\infty]
 h {Q}_\infty^{\pm,1}=
 \frac{\pm h |b_\mp| \sigma_\pm^{2}}{2|b_\pm |(b_-+|b_+|)}.\end{equation}
\end{lemma}
\begin{proof}
Lemma \ref{lemma:ergodic}-{I} implies that 
 \[
 Q^\pm_{h,N}/N \convas[N\to\infty]
 h \EE[ \ind{\{ \pm X_\infty \geq 0\}}] , 
 \]
which, with \eqref{def:Q:stat}, proves \eqref{erg:Q:LLN}. Asymptotic normality \eqref{erg:Q:CLT} also holds because of Lemma \ref{lemma:ergodic}-{II}. 
The limit \eqref{erg:Q1:LLN} follows from Lemma \ref{lemma:ergodic}-{I} and \eqref{def:Q:stat}.
\end{proof}

We are going to prove Lemma~\ref{lemma:lt}. We first need an auxiliary result.
Let us write 
\begin{equation}\label{eq:def:G}
\begin{split}
G(x,h)&:=
\EE[|X_{\infty}^h| \ind{\{ X_{\infty}X_{\infty}^h<0\}} \,\big|X_\infty=x]\\
J(x,h)&:=
\EE[\ind{\{ X_{\infty}X_{\infty}^h<0\}} \,\big|X_\infty=x]
=
\PP(X_{\infty}X_{\infty}^h<0 \,\big|X_\infty=x).
\end{split}
\end{equation}
We have the following lemma.
\begin{lemma}\label{th:laplace}
The Laplace transform
$\mathcal{L}_\lambda G(x,\cdot)=\int_0^\infty e^{-\lambda t} G(x,t)dt $ 
is
\[
\begin{split}
\mathcal{L}_\lambda G(x,\cdot)=
\begin{cases}
\frac{2\sigma_-^2}{
(b_- +\sqrt{b_-^2+2\sigma_-^2\lambda})^2
}
\frac{\exp\left(x\frac{-b_+-\sqrt{b_+^2+2\sigma_+^2 \lambda}}{\sigma_+^2} \right) }{\frac{-b_-+\sqrt{b_-^2+2\sigma_-^2\lambda}}{\sigma_-^2}+
\frac{b_+ +\sqrt{b_+^2+2\sigma_+^2\lambda}}{\sigma_+^2}}
\quad \mbox{ for } x>0 \\
 \frac{2 \sigma_+^2}{
(|b_+|+\sqrt{b_+^2+2\sigma_+^2\lambda})^2 }
 \frac{
 \exp\left(x\frac{-b_-+\sqrt{b_-^2+2\sigma_-^2\lambda}}{\sigma_-^2}\right)
 }
{
\frac{-b_-+\sqrt{b_-^2+2\sigma_-^2\lambda}}{\sigma_-^2}+
\frac{b_+ +\sqrt{b_+^2+2\sigma_+^2\lambda}}{\sigma_+^2}
}
\quad \mbox{ for } x<0
\end{cases}
\end{split}
\]
and the Laplace transform
$\mathcal{L}_\lambda J(x,\cdot)=\int_0^\infty e^{-\lambda t} J(x,t)dt $ 
is
\[
\begin{split}
\mathcal{L}_\lambda J(x,\cdot)=
\begin{cases}
\frac{2}{
b_- +\sqrt{b_-^2+2\sigma_-^2\lambda}
}
\frac{\exp\left(x\frac{-b_+-\sqrt{b_+^2+2\sigma_+^2 \lambda}}{\sigma_+^2} \right) }{\frac{-b_-+\sqrt{b_-^2+2\sigma_-^2\lambda}}{\sigma_-^2}+
\frac{b_+ +\sqrt{b_+^2+2\sigma_+^2\lambda}}{\sigma_+^2}}
\quad \mbox{ for } x>0 \\
 \frac{2}{
|b_+|+\sqrt{b_+^2+2\sigma_+^2\lambda} }
 \frac{
 \exp\left(x\frac{-b_-+\sqrt{b_-^2+2\sigma_-^2\lambda}}{\sigma_-^2}\right)
 }
{
\frac{-b_-+\sqrt{b_-^2+2\sigma_-^2\lambda}}{\sigma_-^2}+
\frac{b_+ +\sqrt{b_+^2+2\sigma_+^2\lambda}}{\sigma_+^2}
}
\quad \mbox{ for } x<0.
\end{cases}
\end{split}
\]
\end{lemma}
\begin{proof}
We apply \eqref{eq:propr:res}, with $f(x,y)=\ind{\{ xy<0\}}|y|$. It is possible to exchange the order of integration because $f(x,y)\geq 0$. Using \eqref{resolventmin} we get, for $x<0$
\begin{equation}
\int_\RR r(\lambda,x,y)f(x,y)dy = 
\frac{\psi_\lambda(x)}{W_\lambda}
\int_0^\infty m(y)\varphi_\lambda(y) f(x,y)
dy .
\end{equation}
The last integral, using \eqref{inv_dens} and \eqref{decreasingsolution}, can be evaluated as
\[
\dfrac{2}{\sigma_+^2}
\int_0^\infty 
	\exp\left(y\frac{-|b_+|-\sqrt{b_+^2+2\sigma_+^2\lambda}}{\sigma_+^2}\right)	
 ydy 
 = 
 \frac{2 \sigma_+^2}{
(|b_+|+\sqrt{b_+^2+2\sigma_+^2\lambda})^2 
 }
\]
so that
\[
\cL_\lambda \EE_x [f(x,X_\cdot)]
=
 \frac{2 \sigma_+^2}{
(|b_+|+\sqrt{b_+^2+2\sigma_+^2\lambda})^2 }
 \frac{
 \exp\left(x\frac{-b_-+\sqrt{b_-^2+2\sigma_-^2\lambda}}{\sigma_-^2}\right)
 }
{
\frac{-b_-+\sqrt{b_-^2+2\sigma_-^2\lambda}}{\sigma_-^2}+
\frac{b_+ +\sqrt{b_+^2+2\sigma_+^2\lambda}}{\sigma_+^2}
}.
\]
For $x>0$ the computation is analogous. If we now apply \eqref{eq:propr:res}, with $f(x,y)= \ind{\{xy<0\}}
$, which again is positive, using \eqref{resolventmin} we get, for $x<0$
\begin{equation}
\int_\RR r(\lambda,x,y)f(x,y)dy = 
\frac{\psi_\lambda(x)}{W_\lambda}
\int_0^\infty m(y)\varphi_\lambda(y) 
dy.
\end{equation}
The last integral, using \eqref{inv_dens} and \eqref{decreasingsolution}, can be evaluated as
\[
\dfrac{2}{\sigma_+^2}
\int_0^\infty 
	\exp\left(y\frac{-|b_+|-\sqrt{b_+^2+2\sigma_+^2\lambda}}{\sigma_+^2}\right)	
 dy 
 = 
 \frac{2 }{
|b_+|+\sqrt{b_+^2+2\sigma_+^2\lambda} 
 }
\]
so that
\[
\cL_\lambda \EE_x [f(x,X_\cdot)]
=
 \frac{2}{
|b_+|+\sqrt{b_+^2+2\sigma_+^2\lambda} }
 \frac{
 \exp\left(x\frac{-b_-+\sqrt{b_-^2+2\sigma_-^2\lambda}}{\sigma_-^2}\right)
 }
{
\frac{-b_-+\sqrt{b_-^2+2\sigma_-^2\lambda}}{\sigma_-^2}+
\frac{b_+ +\sqrt{b_+^2+2\sigma_+^2\lambda}}{\sigma_+^2}
}.
\]
For $x>0$ the computation is analogous.
\end{proof}

\subsubsection{Proof of Lemma \ref{lemma:lt}}
Let us start with $L_{h,N}$. Writing $g(x,y)=\ind{\{ xy<0\}}|y|$, from \eqref{def:LT:discr} we have
$
	L_{h,N}= 
2
	 \sum_{k=0}^{N-1} g( \X_k^h).
$
Lemma \ref{lemma:ergodic}-{III} implies
 that 
 \[
 L_{h,N}/N \convas[N\to\infty]
2 \EE[g(\X_\infty^h)],
 \]
so that, using the tower property of conditional expectation, we get
\[
	\EE[g(\X_\infty^h)]  
	= \EE[|X_{\infty}^h| \ind{\{ X_{\infty}X_{\infty}^h<0\}}]
	=\EE[G(X_\infty,h)]
\]
with $G$ in \eqref{eq:def:G}. 
The proof strategy consists in computing this quantity via Laplace transform.
The Laplace transform 
$
\mathcal{L}_\lambda G(x,\cdot)=\int_0^\infty e^{-\lambda t} G(x,t)dt 
$ 
is computed in Lemma \ref{th:laplace}. Using the stationary distribution \eqref{inv_dens}, explicit computations give
\[
\begin{split}
\EE[\mathcal{L}_\lambda G(X_\infty,\cdot)] 
&=\frac{1}{M(\RR)}\int_{-\infty}^0 \mathcal{L}_\lambda G(x,\cdot) m(x) dx
+\frac{1}{M(\RR)}\int_0^\infty \mathcal{L}_\lambda G(x,\cdot) m(x) dx
\\
&
= \dfrac{|b_+| b_-}{b_-+|b_+|}
\frac{1
}{\lambda^2}.
\end{split}
\]
By Tonelli's theorem, it holds that
\[
\EE[\mathcal{L}_\lambda G(X_\infty,\cdot)] = \int_0^\infty e^{-\lambda t} \EE[G(X_\infty,t)]dt = \mathcal{L}_\lambda \EE[ G(X_\infty,\cdot)].
\] 
Therefore, inverting the Laplace transform, we obtain
\begin{equation}\label{eq:Eg}
\EE[g(\X_\infty^h)]=
\EE[G(X_\infty,h)]=
\mathcal{L}^{-1}_h\left(
\dfrac{|b_+| b_-}{b_-+|b_+|}
\frac{1
}{(\cdot)^2}\right)
=
\dfrac{|b_+| b_-}{b_-+|b_+|} h,
\end{equation}
which proves \eqref{eq:cons:lt}.

Let us now consider $\hat{L}_{h,N}$. One can rewrite $\hat{L}_{h,N}$ as
\[
\hat{L}_{h,N}=	
-\frac{3\sqrt{\pi}}{2\sqrt{2}}
\frac{\sigma_-+\sigma_+}{\sigma_-\sigma_+}
  \frac{1}{\sqrt{h}}\sum_{k=0}^{N-1} X_{kh
	 } X_{(k+1)h}\ind{\{X_{kh
	 } X_{(k+1)h}<0\}}
\]
so that Lemma \ref{lemma:ergodic}-{III} implies that as $N\to\infty$, ${\hat{ L}_{h,N}}/{N}$ converges $\PP$-a.s.~to
 \[
-\frac{3\sqrt{\pi}}{2\sqrt{2}}
\frac{\sigma_-+\sigma_+}{\sigma_-\sigma_+}
  \frac{1}{\sqrt{h}}
  \EE[|X_{\infty}| g(\X_\infty^{h})]
  =
-\frac{3\sqrt{\pi}}{2\sqrt{2}}
\frac{\sigma_-+\sigma_+}{\sigma_-\sigma_+}
 \EE[X_{\infty}X_{\infty}^h \ind{\{ X_{\infty}X_{\infty}^h<0\}}]. 
 \]
Using Lemma \ref{th:laplace} and the stationary distribution \eqref{inv_dens}, explicit computations give
\[
\begin{split}
& \EE[\mathcal{L}_\lambda  (|X_\infty| G(X_\infty,\cdot))] 
=\frac{1}{M(\RR)}\int_{-\infty}^0 \mathcal{L}_\lambda G(x,\cdot) m(x) |x| dx
+\frac{1}{M(\RR)}\int_0^\infty \mathcal{L}_\lambda G(x,\cdot) m(x) |x| dx
\\&
=
 \frac{|b_{+}|b_{-}}{|b_{+}|+b_{-}} \frac{1}{\lambda^{3}}
\frac{(-b_-+\sqrt{b_-^2+2\sigma_-^2\lambda})(b_+ + \sqrt{b_+^2+2\sigma_+^2\lambda})}{
b_-+\sqrt{b_-^2+2\sigma_-^2\lambda}+|b_+| + \sqrt{b_+^2+2\sigma_+^2\lambda}
}
\\&
=\sqrt{2} \frac{|b_{+}|b_{-}}{|b_{+}|+b_{-}}
\frac{\sigma_-\sigma_+}{
\sigma_-+\sigma_+}
 \frac{1}{\lambda^{5/2}}
 +
  \frac{|b_{+}|b_{-}}{|b_{+}|+b_{-}}
 \left( 
  \frac{\sigma_-b_{+}-\sigma_+b_{-}}{
\sigma_-+\sigma_+}
  -
\frac{(b_{-} +|b_{+}| ) \sigma_-\sigma_+}{
(\sigma_-+\sigma_+)^{2}}
  \right) 
 \frac{1}{\lambda^{3}}
+
O( \lambda^{-7/2}) 
\end{split}
\]
where with $O(\lambda^{-7/2})$ we denote a function $f(\lambda)$ s.t. $\lambda^{7/2} f(\lambda)\to C\in \RR,$ as $ \lambda \to \infty$.
Therefore, inverting the Laplace transform, using linearity and  a Tauberian theorem for $O(\lambda^{-7/2})$ we obtain
\begin{equation}\label{eq:Eg}
\begin{split}
&\EE[|X_\infty| g(\X_\infty^h)]=
\mathcal{L}^{-1}_h\EE[\mathcal{L}_\lambda (|X_\infty| G(X_\infty,\cdot))] 
\\
&=
 \frac{\sqrt{2}}{\Gamma(5/2)}
\frac{|b_{+}|b_{-}}{|b_{+}|+b_{-}}
\frac{\sigma_-\sigma_+}{
\sigma_-+\sigma_+}
 h^{3/2}
 \\& 
 +
  \frac{1}{\Gamma(3)} 
  \frac{|b_{+}|b_{-}}{|b_{+}|+b_{-}}
 \left( 
  \frac{\sigma_-b_{+}-\sigma_+b_{-}}{
\sigma_-+\sigma_+}
  -
\frac{(b_{-} +|b_{+}| ) \sigma_-\sigma_+}{
(\sigma_-+\sigma_+)^{2}}
  \right) 
  h^{2}
 +O(h^{{5/2}})
 \\&
=: 
\hat{L}(h) \left(
 \frac{2\sqrt{2}}{3 \sqrt{\pi}}
\frac{\sigma_-\sigma_+}{
\sigma_-+\sigma_+}
h^{{3/2}} \right)
 \end{split}
 \end{equation}
as $h \to 0$, where $\Gamma$ denotes the $\Gamma$-function. This implies~\eqref{eq:cons:lt-2}.

Let us now consider $\bar{L}_{h,N}$. Let $j(x,y):= \ind{\{xy<0\}}$, then 
$\bar{L}_{h,N} = \sqrt{\frac{\pi}2} \frac{\sigma_++\sigma_-}{2} \sqrt{h} \sum_{k=0}^{N-1} j( \X_k^h)$.
The proof works the same as in the case of $L_{h,N}$ but $g,G$ are now $j,J$.
Lemma~\ref{th:laplace}, the stationary distribution \eqref{inv_dens} and explicit computations give
\[
\begin{split}
\EE[\mathcal{L}_\lambda J(X_\infty,\cdot)] 
&=\frac{1}{M(\RR)}\int_{-\infty}^0 \mathcal{L}_\lambda J(x,\cdot) m(x) dx
+\frac{1}{M(\RR)}\int_0^\infty \mathcal{L}_\lambda J(x,\cdot) m(x) dx
\\
&
= \dfrac{8 |b_+| b_-}{b_-+|b_+|}
\frac{1}{2\lambda}
\frac{1}{|b_+| +\sqrt{b_+^2+2\sigma_+^2\lambda}
+b_-+\sqrt{b_-^2+2\sigma_-^2\lambda}}
 \\
 &
= \sqrt{2} \dfrac{ 2|b_+| b_-}{b_-+|b_+|}
\frac{1}{\sigma_+ + \sigma_-}
 \frac{1}{\lambda^{{3/2}}}
 -
 \dfrac{ 2|b_+| b_-}{(\sigma_+ + \sigma_-)^{2}}
 \frac{1}{\lambda^2}
+O(\lambda^{-5/2}) 
\end{split}
\]
as $\lambda \to \infty$.
Therefore, inverting the Laplace transform and using a Tauberian theorem like before we obtain
\begin{equation}\label{eq:Eg}
\begin{split}
\EE[j(\X_\infty^h)]& = \mathcal{L}^{-1}_h\EE[\mathcal{L}_\lambda J(X_\infty,\cdot)] 
\\
&=\frac{1}{
\Gamma(3/2)} \dfrac{ |b_+| b_-}{b_-+|b_+|}
\frac{2\sqrt{2}}{\sigma_++\sigma_-}
\sqrt{h}
+
\frac{1}{
\Gamma(2)}
\dfrac{ 2|b_+| b_-}{(\sigma_+ + \sigma_-)^{2}}
h
+O(h^{3/2})
\\
&=: 
\bar{L}(h)\left(\frac{\sqrt{2}}{
\sqrt{\pi}} \frac{2}{\sigma_++\sigma_-}
 h^{{1/2}}\right)
\end{split}
\end{equation}
which implies \eqref{eq:cons:lt-1}.

\subsubsection{Proof of Theorem~\ref{th:fixed:h} and Theorem~\ref{th:sigma:h}}
\begin{lemma}\label{lem:jointconv}
It holds
\[
 \frac{1}{\sqrt{N}}
 \begin{pmatrix}
 \frac{1}{ h Q_\infty^+}\left(
 -L_{h,N}/2
- Q_{h,N}^+ b_+ \right)\\
\frac{1}{ h Q_\infty^-}
\left(L_{h,N}/2
- Q_{h,N}^- b_-\right)
\end{pmatrix} 
\convl[N\to\infty]
 \cN(0,\Gamma),
\]
where
\begin{equation}\label{eq:def:gammaq}
\begin{split}
\Gamma=\EE[q(\tilde{\X}^{h}_0) q^{T}(\tilde{\X}^{h}_0) ]+\sum_{k=1}^\infty  \EE[q(\tilde{\X}^{h}_0)q^{T}(\tilde{\X}^{h}_k)]
+
\EE[q(\tilde{\X}^{h}_k)q^{T}(\tilde{\X}^{h}_0)]
\end{split}
\end{equation}
with
\begin{equation}\label{eq:q}
q(x,y)=  
\begin{pmatrix}
q_{1}(x,y)\\
q_{2}(x,y)\end{pmatrix} 
=
\begin{pmatrix}
\frac{b_-+|b_+|}{b_-}
(
-\ind{\{ xy<0\}} |y|
- b_+ h \ind{\{x>0\}})\\
\frac{b_-+|b_+|}{b_+}
(
\ind{\{ xy<0\}} |y|
- b_- h \ind{\{x<0\}}
)
\end{pmatrix} .
\end{equation}
\end{lemma}
\begin{proof} Note that
\[
 \begin{pmatrix}
 \frac{1}{ h Q_\infty^+}\left(
 -L_{h,N}/2
- Q_{h,N}^+ b_+ \right)\\
\frac{1}{ h Q_\infty^-}
\left(L_{h,N}/2
- Q_{h,N}^- b_-\right)
\end{pmatrix} 
=
\sum_{k=0}^{N-1}
 q( \X_k^h).
\]
 From
\eqref{eq:Eg} and \eqref{def:Q:stat} we have $\EE[q(\X_\infty^h)]=0$. So, Lemma \ref{lemma:ergodic}-{III} gives
\[
\frac{1}{N}
\sum_{k=0}^{N-1}
 q( \X_k^h)\convas[N\to\infty] 0, 
\]
and Lemma \ref{lemma:ergodic}-{IV}, using again $\EE[q(\X_\infty^h)]=0$, gives the convergence in the statement.
\end{proof}

\begin{proof}[Proof of Theorem \ref{th:fixed:h}] 
Consistency is clear from \eqref{erg:Q:LLN} and \eqref{eq:cons:lt}.
The limit error of estimator $(\overline{b}^-_{h,N},\overline{b}^+_{h,N})$ can be rewritten with some manipulations as
\[
\sqrt{N}
\begin{pmatrix}
 \overline{b}^+_{h,N}
-b_+ \\
 \overline{b}^-_{h,N}
-b_-
\end{pmatrix}
=
 \frac{1}{\sqrt{N}}
\begin{pmatrix}
\frac{1}{ h Q_\infty^+}\left(
 -\frac{L_{h,N}}{2 }
- Q_{h,N}^+b_+\right)
\\
\frac{1}{ h Q_\infty^-} \left(
 \frac{L_{h,N}}{2 }
- Q_{h,N}^-b_-\right)
\end{pmatrix}
+
\begin{pmatrix}
\frac{\sqrt{N}}{2hQ_\infty^+}
\left(hQ_\infty^+ - \frac{Q_{h,N}^+}{N} \right)
 \left(
- \frac{L_{h,N}}{Q_{h,N}^+}
+
 \frac{ \Linfty}{  Q_\infty^+}
\right)
\\
\frac{\sqrt{N}}{2hQ_\infty^-}
\left(hQ_\infty^- - \frac{Q_{h,N}^-}{N} \right)
 \left(
 \frac{L_{h,N}}{Q_{h,N}^-}
-
 \frac{ \Linfty}{  Q_\infty^-}
\right)
\end{pmatrix}.
\]
The last summand vanishes because of \eqref{erg:Q:LLN}, \eqref{erg:Q:CLT}, and \eqref{eq:cons:lt}.
The convergence of the first summand follows from Lemma \ref{lem:jointconv} and implies Theorem \ref{th:fixed:h}-ii)  (also using \eqref{def:Q:stat}). 
\end{proof}

\begin{proof}[Proof of Theorem~\ref{th:sigma:h}] 
Consistency follows from \eqref{erg:Q:LLN}, \eqref{erg:Q1:LLN} and \eqref{eq:cons:lt}. To prove asymptotic normality of the positive side we write
\[
\begin{split}
\sqrt{N}
 ( (& \bar{\sigma}^+_{h,N})^2
-\sigma_+^2) 
=
\frac{\sqrt{N}}{2} 
\left(
\frac{L_{h,N} Q^{+,1}_{h,N}}{ ({Q}^{+}_{h,N})^2} 
-
 \frac{L_{\infty} Q^{+,1}_{\infty}}{ ({Q}^{+}_{\infty})^2} 
\right)\\
&=\frac{\sqrt{N}}{2} 
\frac{L_{h,N}/N Q^{+,1}_{h,N}/N (h{Q}^{+}_{\infty})^2-hL_{\infty} hQ^{+,1}_{\infty}({Q}^{+}_{h,N}/N)^2}
{ ({Q}^{+}_{h,N}/N \,h{Q}^{+}_{\infty})^2} 
\\
&
=\frac{\sqrt{N}}{2} 
\frac{
(L_{h,N}/N-h \Linfty)  Q^{+,1}_\infty( {Q}^{+}_{\infty})^2
+
 \Linfty (Q^{+,1}_{h,N}/N - h Q^{+,1}_{\infty})( {Q}^{+}_{\infty})^2
}{ h ( {Q}^{+}_{\infty})^4} 
 + \\
&
\quad\quad+\frac{\sqrt{N}}{2} 
\frac{
2  \Linfty Q^{+,1}_{\infty}
{Q}^{+}_{\infty}
(h {Q}^{+}_{\infty}
- {Q}^{+}_{h,N}/N)
}{ h ( {Q}^{+}_{\infty})^4} 
 + o_N(1)
  \\
&
=\frac{\sqrt{N}}{2} 
\left(
\frac{L_{h,N}}{N}  
\frac{\sigma_+^2(b_-+|b_+|)}{2h |b_+| b_-}
+
\frac{Q^{+,1}_{h,N}}{N}
\frac{2|b_+| (b_-+|b_+|)}{hb_-}
- \frac{{Q}^{+}_{h,N}}{N}
\frac{2\sigma_+^2(b_-+|b_+|)}{h  b_-} \right)
 + o_N(1)
 \end{split}
 \] 
where with $o_N(1)$ we mean a quantity that goes to $0$ in probability as $N\to\infty$, that may vary from line to line. A similar computation also holds for the negative side, so that we can write
\[
\begin{split}
&\sqrt{N}
\begin{pmatrix}
 \bar{\sigma}^+_{h,N}
-\sigma^{2}_+ \\
  \bar{\sigma}^-_{h,N}
-\sigma^{2}_-
\end{pmatrix}
=\frac{1}{\sqrt{N}}
\sum_{k=0}^{N-1}
v(\X_k^h) 
 + o_N(1)
\end{split}
\]
with
\[
v(x,y)
=
\frac{(b_-+|b_+|)}{b_-|b_+|} 
\left(
\ind{\{xy<0\}}
\frac{|y|}{2h} 
\begin{pmatrix}
\sigma_+^2\\
\sigma_-^2
\end{pmatrix}
+
|x|
\begin{pmatrix}
b_+^2 \ind{\{x>0\}}\\
b_-^2 \ind{\{x<0\}}
\end{pmatrix}
- 
\begin{pmatrix}
\sigma_+^2 |b_{+}| \ind{\{x>0\}}\\
\sigma_-^2 b_{-} \ind{\{x<0\}}
\end{pmatrix}
\right).
\]
It can be verified from \eqref{eq:Eg} and \eqref{def:Q:stat} that $\EE[v(\X_\infty^h) ]=0$. Using Lemma \ref{lemma:ergodic}-IV allows to conclude, analogously to the proof of Lemma \ref{lem:jointconv}, this time with asymptotic covariance
\begin{equation}\label{eq:def:gammav}
\begin{split}
\Gamma_{\sigma}=\EE[v(\tilde{\X}^{h}_0) v^{T}(\tilde{\X}^{h}_0) ]+\sum_{k=1}^\infty  \EE[v(\tilde{\X}^{h}_0)v^{T}(\tilde{\X}^{h}_k)]
+
\EE[=v(\tilde{\X}^{h}_k)v^{T}(\tilde{\X}^{h}_0)],
\end{split}
\end{equation}
so that the theorem is proved.
\end{proof}

\subsection{Proof of Theorem~\ref{th:joint:CLT}: maximum likelihood estimator, high frequency and long time horizon}\label{sec:mle}


In \cite{lp2} it was proved that the MLE from continuous time observation on the time interval $[0,T]$, $T\in (0,\infty)$ is given by~
\begin{equation} \label{eq:QMLE_ct}	
	\beta_T:=\begin{pmatrix}
	\beta^+_T, & \beta^-_T 
	\end{pmatrix}
	= 
	\begin{pmatrix}
	\frac{	\cM^{+}_T }{ 
	Q^{+}_T }
	, & \frac{	\cM^{-}_T }{ 
	Q^{-}_T } 
	\end{pmatrix}.
\end{equation}
Note that it coincides with the estimator obtained my maximizing the quasi-likelihood function, which is defined as the likelihood, but with diffusion coefficient set to $1$.
The discrete version of the estimator has been provided as well and is~\eqref{eq:QMLE_dis_time} as the quantity maximizing the \emph{discretized likelihood}:
\begin{equation*}
    G_{h,N}(b_+,b_-)=\exp\bigg(\sum_{i=0}^{N-1} \frac{b(X_{i h})}{(\sigma(X_{i h}))^2}(X_{(i+1)h}-X_{i h})-\frac{h}{2}\sum_{i=0}^{N-1} \frac{(b(X_{i h}))^2}{(\sigma(X_{ i h}))^2} \bigg).
\end{equation*}

We now state the results about the asymptotic behavior in long time of the estimators based on continuous observations.
\begin{lemma}[cf.~Proposition~3 in~\cite{lp2}] \label{th:continuous}	
Assume the process is ergodic.
Let $\pm\in \{+,-\}$. The MLE~\eqref{eq:QMLE_ct}
\begin{enumerate}[i)]
\item \label{th1:item:LLN}
is consistent:
$ 
\beta^\pm_T-b_\pm \convas[T\to\infty] 0,
$
and
\item \label{th1:item:CLT}
satisfies the following CLT: 
$
\sqrt{T}
( \beta^\pm_T-b_\pm )
\xrightarrow[{T\to \infty}]{\mathrm{stably}}  
N^\pm
$
where  $N^+,N^-$ are two mutually independent, independent of $X$, Gaussian random variables 
with variance respectively $\sigma_\pm^2 \frac{|b_+| +b_-}{|b_\mp|}$.
\end{enumerate}
\end{lemma}


The proof of Theorem~\ref{th:joint:CLT} exploits the latter lemma and the key Lemma~\ref{lemma:control:diff:disc:cont} which is based on a bound for the transition density proposed in Lemma~\ref{lem:skew:density}.
More precisely for all $N\in \NN$ it holds
\[
 \widehat{b}^\pm_{h_N,N}
-b_\pm
=
 \widehat{b}^\pm_{h_N,N}
-\beta_{T_N}^\pm
+
\beta_{T_N}^\pm-b_\pm.
\]
We then handle the second term of the sum with Lemma~\ref{th:continuous}. 
Equations~\eqref{eq:QMLE_ct} and~\eqref{eq:QMLE_dis_time}
imply for $j\leq 1$ that 
\begin{equation} \label{eq:difference:CLT}
\begin{split}
&	\sqrt{T_N^{j}} \left( \widehat{b}^\pm_{h_N,N} -\beta^\pm_{T_N}\right) 
 = \frac{	\cM^{\pm}_{h_N,N} - \cM^{\pm}_{T_N} }{ \sqrt{T_N} }
\frac{\sqrt{T_N^{1+j}}}{ 
\cQ^{\pm}_{h_N,N} } 
+ \frac{\cM^{\pm}_{T_N}}{\sqrt{T_N}} \sqrt{T_N^{1+j}}\left( \frac{1}{ 
\cQ^{\pm}_{h_N,N} }- \frac{1}{ 
\cQ^{\pm}_{T_N} }\right)
\\
& =
\frac{	\cM^{\pm}_{h_N,N} - \cM^{\pm}_{T_N} }{ \sqrt{T_N} }
\frac{\sqrt{T_N^{1+j}}}{ 
\cQ^{\pm}_{h_N,N} } 
+ \frac{\cM^{\pm}_{T_N}}{\sqrt{T_N^{3-j}}}  \frac{ \cQ^{\pm}_{T_N}-\cQ^{\pm}_{h_N,N}}{ \sqrt{T_N} } \frac{T_N^2 }{ 
\cQ^{\pm}_{h_N,N} \cQ^{\pm}_{T_N} }.
\end{split}
\end{equation}
The latter converges in probability to 0 by Lemma~\ref{lemma:control:diff:disc:cont} and the fact that ${\cQ^{\pm}_{T_N}}/{T_N}$ converges in probability (ergodicity) and ${\cM^{\pm}_{T_N}}/{\sqrt{T_N}} $ converges in law (ergodicity and martingale CLT in~\cite{crimaldi}).
Lemma~\ref{lemma:control:diff:disc:cont} is similar to Lemma~4 in~\cite{mazzonetto2020drift}. Nevertheless, here we allow a deterministic initial condition using the density bound in Lemma~\ref{lem:skew:density}.

\begin{lemma}
\label{lemma:control:diff:disc:cont}
Let $\alpha\in [0,1]$ be fixed and let $X$ be the solution to \eqref{eq:DOBM}, with $X_0$ deterministic, 
and let $(T_N)_{N\in \NN},(h_N)_{N\in \NN}$ be positive sequences satisfying $h_N=T_N/N$ such that, as $N \to\infty$, that $T_N\to \infty$ and 
$
	T_N^{2(1-\alpha)}h_N \to 0.
$
Then
\begin{equation*}
\limsup_{N\to\infty} {T_N^{-\alpha}}
 \EE\left[ | \cQ^{\pm}_{T_N} - \cQ^{\pm}_{h_N,N} |\right] = 0
\quad \text{and} \quad
\limsup_{N\to\infty} {T_N^{-\alpha}} \EE\left[ |\cM^{\pm}_{T_N} - \cM^{\pm}_{h_N,N} | \right] = 0
\end{equation*}
where $\cQ^{\pm}_{T_N}$, $\cQ^{\pm}_{h_N,N}$,
$\cM^{\pm}_{T_N}$, $\cM^{\pm}_{h_N,N}$ are defined in~\eqref{def:M:Q} and \eqref{def:M:Q:disc}.
\end{lemma}

\begin{proof}[Proof of Lemma~\ref{lemma:control:diff:disc:cont}]
In this proof we use the round ground notation 
$\lfloor t \rfloor_{N} := t_k$ if $t \in [t_k , t_{k+1})$, where $t_k := k h_N$.
Moreover, without loss of generality, we assume $1\leq T_N \leq N$ for all $N\in \NN$.

Let us first note for every $t\in (0,\infty)$ that
\[
	\ind{\{ \pm X_{t} >0 \}} - \ind{\{ \pm X_{\lfloor t \rfloor_{N}} >0 \}}
	= \mp \sgn(X_{\lfloor t \rfloor_{N}})  \ind{\{X_{\lfloor t \rfloor_{N}} X_t<0\}}.
\]
Hence
\begin{equation} \label{eq:diffQ_TN}
\begin{split}
 \EE\left[|\cQ^{\pm}_{T_N} - \cQ^{\pm}_{h_N,N}|\right]
	& \leq \int_0^{T_N} \EE\!\left[\ind{\{X_{\lfloor t \rfloor_{N}} X_t<0\}} \right] \vd t
	= \int_0^{T_N} \PP\!\left( X_{\lfloor t \rfloor_{N}} X_t<0 \right) \vd t .
\end{split} 
\end{equation}
Analogously, triangular inequality, H\"older's inequality, and It\^o-isometry imply that
\begin{equation} \label{eq:diffM_TN}
\begin{split}
	& \EE\left[ |\cM^{\pm}_{T_N} - \cM^{\pm}_{h_N,N} | \right]
	\\
	& \leq  (|b_-| \vee |b_+|) \int_0^{T_N} \EE\left[ \ind{\{X_t X_{\lfloor t \rfloor_{N}}< 0\}} \right] \vd t
	+ 
	(\sigma_-\vee \sigma_+) \left(\int_0^{T_N}  \EE\left[  \ind{\{X_t X_{\lfloor t \rfloor_{N}}< 0\}} \right] \vd t\right)^{\!\nicefrac12}.
\end{split}
\end{equation}

Hence, the proof of Lemma~\ref{lemma:control:diff:disc:cont}, which consists in showing that
\eqref{eq:diffQ_TN} and \eqref{eq:diffM_TN} are $o(T_N^\alpha)$, 
reduces to prove 
\begin{equation} \label{item:cond:1}
	\int_0^{T_N} \PP\!\left( X_{\lfloor t \rfloor_{N}} X_t<0 \right)  \!\vd t 
\text{ is } o(T_N^{\alpha}).
\end{equation}

{\it Step 1}:  {\it Let $s,t\in (0,\infty)$ be fixed such that $s\leq t$.
We show that there exist a constant $C \in (0,\infty)$ depending only on $b_\pm,\sigma_\pm$ such that
for $N$ large enough 
\begin{equation} \label{th:joint:step2}
\EE\!\left[   \ind{\{X_{s} X_t<0\}} \right] 
\leq C \left(1+\frac1{\sqrt{s}}\right) \sqrt{t-s}.
\end{equation}
}

In order to prove this, let us fix $\pm\in \{-,+\}$. Note that $\EE\left[\EE\left[   \ind{\{ \pm X_t<0\}} \ind{\{\pm X_{s} > 0\}}  |  X_{s} \right]\right]$
is bounded by
$ \EE\left[\PP\left( \tau_s^\pm \leq t-s \right) \right]$ 
where $\tau_s^\pm$ is the first hitting time of the level 0 of the process solution to the Brownian motion with drift
\[ \xi_u = X_s + b_\pm u +  \sigma_\pm  W_u^s = X_s + b(X_s) u +  \sigma(X_s)  W_u^s \]
with $W^s$ a Brownian motion independent of $\sigma(X_v, v\in [0,s])$.
Then 
\[
\begin{split}
	\EE\left[   \ind{\{ \pm X_t<0\}} \ind{\{\pm X_{s} > 0\}}  |  X_{s} \right] 
\leq \PP\left( \tau_s^\pm \leq t-s \right)  
&= 
\int_0^{t-s} \frac{|X_s|}{\sigma_\pm\sqrt{2\pi u^3}} 
\exp\left(
	-\frac{(X_s-b_\pm u)^2}{2\sigma_\pm^2 u}
\right)	
\vd u.
\end{split}
\]
Let $f_\pm(y;u) := \exp\left(	-\frac{y^2}{2\sigma_\pm^2 u}\right)$.
We apply Lemma~\ref{lem:skew:density} and obtain for some non-negative constant $C$:
\[
\begin{split}
	& \EE\left[   \ind{\{ \pm X_t<0\}} \ind{\{\pm X_{s} > 0\}} \right] 
	\leq \int_{\RR_\pm}  q_{(b_-,b_+)}(s,X_0,x) \int_0^{t-s} \frac{|x|}{\sigma_\pm\sqrt{2\pi u^3}} f_\pm(x-b_\pm u ;u) \vd u  \vd x
	\\ &
	\leq C \frac{1+\sqrt{s}}{\sqrt{s}}    \int_0^{t-s} \int_{\RR_\pm} \frac{|x+2b_\pm u|}{\sigma_\pm^2\sqrt{2\pi u^3}} f_\pm(x;u) \vd x \vd u 
	\\
	& \leq  C \frac{1+\sqrt{s}}{\sqrt{s}}   \int_0^{t-s} \int_{\RR_\pm} \frac{|x|}{\sqrt{2\pi u^3}} e^{-\frac{x^2}{2u}} \vd x \vd u 
	+ \frac{2 |b_\pm|}{\sigma_\pm^2} C \frac{1+\sqrt{s}}{\sqrt{s}}   \int_0^{t-s} \int_{\RR_\pm} \frac{1}{\sqrt{2\pi u}} e^{-\frac{x^2}{2u}} \vd x \vd u
	\\ &
	= C \frac{1+\sqrt{s}}{\sqrt{s}} \left(\frac1{\sqrt{2\pi}}  \int_0^{t-s} \frac1{\sqrt{u}} \vd u + \frac{2 |b_\pm|}{\sigma_\pm^2}  (t-s) \right)
	\\ & 
	= \left( \sqrt{\frac{2}{\pi}} + \frac{2 |b_\pm|}{\sigma_\pm^2} \sqrt{t-s}\right) C \left(1+\frac1{\sqrt{s}}\right) \sqrt{t-s}.
\end{split}
\]
This completes the proof of the first step.

{\it Step 2}: {\it Proof of~\eqref{item:cond:1}}.

For all $t\in [t_1,\infty)$ the tower property and~\eqref{th:joint:step2} imply, in the notation of~\eqref{th:joint:step2}, that 
\[
\begin{split}
	&  \PP\!\left( X_{\lfloor t \rfloor_{N}} X_t<0 \right) 
	=  \EE\!\left[  \EE\!\left[  \ind{\{X_{\lfloor t \rfloor_{N}} X_t<0\}} | X_{\lfloor t \rfloor_{N}} \right] \right]
	\leq
	C \left(1+\frac1{\sqrt{\lfloor t \rfloor_{N}}}\right) \sqrt{t-\lfloor t \rfloor_{N}}.
\end{split}
\]
Then, 
\[
\begin{split}
	\frac{1}{T_N^{\alpha}} & \int_0^{T_N} \PP\!\left( X_{\lfloor t \rfloor_{N}} X_t<0 \right)   \! \vd t
 		\leq 
		\frac{2 h_N}{T_N^{\alpha}} 
		+ C \frac{1}{T_N^{\alpha}} \int_{2h_N}^{T_N}  \sqrt{t-\lfloor t \rfloor_{N}} \! \vd t 
		+ C \frac{1}{T_N^{\alpha}} \int_{2h_N}^{T_N}  \frac{\sqrt{t-\lfloor t \rfloor_{N}}}{\sqrt{\lfloor t \rfloor_{N}}} \! \vd t.
\end{split}
\]
On the right-hand-side of the latter inequality, the first term vanishes as $N\to\infty$, the second term is bounded by 
\[
	C \frac{1}{T_N^{\alpha}} \int_0^{T_N}  \sqrt{t-\lfloor t \rfloor_{N}} \! \vd t 
	\leq C \sqrt{h_N T_N^{2(1-\alpha)}},
\]
and the last term is bounded by
\[
	C \frac{1}{T_N^{\alpha}} \int_{2h_N}^{T_N}  \frac{\sqrt{t-\lfloor t \rfloor_{N}}}{ \sqrt{t-h_N}} \! \vd t 
	\leq C \frac{\sqrt{h_N}}{T_N^{\alpha}} \int_{0}^{T_N}  \frac{1}{ \sqrt{u}} \! \vd u 
	= C \sqrt{h_N T_N^{1-2\alpha}}.
\]
The proof is thus completed.
\end{proof}


Let $q_{(b_+,b_-)}(t,x,y)$ denote the transition density of the drifted oscillating Brownian motion~\eqref{eq:DOBM} with diffusion coefficient $\sigma$ in \eqref{sigmaDOBM}. 
To bound the density, we take inspiration from the article~\cite{downes2009bounds}. 

Note that the transition density of the oscillating Brownian motion without drift, $q_{(0,0)}(t,x,y)$, is known. It was provided in~\cite{kw}, 
and it satisfies for instance that $q_{(0,0)}(t,x,y) \leq C t^{-1/2}$ for some constant $C\in (0,\infty)$ depending on the parameters.

\begin{lemma} \label{lem:skew:density}
There exist $c,C,K \in (0,\infty)$ depending on the parameters $b_\pm,\sigma_\pm$ such that for all $y \in \RR$, $t\in (0,\infty)$
it holds that
\[ q_{(b_+,b_-)}(t,x,y) \leq C e^{y b(y)/\sigma^2(y)- x b(x)/\sigma^2(x)} (1+\sqrt{t}) q_{(0,0)}(t,x,y).\]
Moreover, 
$q_{(0,0)}(t,x,y) \leq \frac1{\sigma(y)} (1+\frac{|\sigma_--\sigma_+|}{\sigma_-+
\sigma_+}) \frac{1}{\sqrt{2\pi t}} e^{-{\left(\frac{y}{\sigma(y)^2}- \frac{x}{\sigma(x)^2}\right)}/{2 t}} .$
\end{lemma}
\begin{proof}[Proof of Lemma~\ref{lem:skew:density}]
Note that, if $Y$ is a drifted oscillating Brownian motion~\eqref{eq:DOBM} with diffusion coefficient $\sigma$~in \eqref{sigmaDOBM} and drift $b(x)$, then $X/\sigma(X)$ is a skew Brownian motion (SBM) with drift $b(x)/\sigma(x)$ and skewness parameter $\frac{\sigma_--\sigma_+}{\sigma_-+\sigma_+}$. This follows from It\^o-Tanaka formula, for a proof see e.g.~\cite{mazzonetto2019rates}.
Moreover, if $p_{(b_+,b_-)}(t,x,y)$ denotes the transition density of the $\beta:=\frac{(\sigma_--\sigma_+)}{(\sigma_-+\sigma_+)}$-SBM with drift $b$, then  
\[q_{(b_+,b_-)}(t,x,y) = p_{(b_+/\sigma_+,b_-/\sigma_-)} (t, x/\sigma(x), y/\sigma(y)) / {\sigma(y)}.\]
In this proof we only consider $X$ a skew-BM with piecewise constant drift $b$ (abuse of notation because $b$ is actually $b/\sigma$) and skewness parameter $\beta$ such that $X_0=x_0$.
By Girsanov's theorem there exists a measure $\QR$ absolutely continuous with respect to $\PP$ such that $X$ under $\QR$ is a driftless $\beta$-SBM with $X_0=x_0$ and for all $f$ measurable and bounded it holds that
\[ \EE\!\left[ f(X_t) \right] 
= \EE_{\QR}\!\left[  \exp{\!\left( \int_0^t  b(X_s) \vd W_s - \frac12 \int_0^t  (b(X_s))^2 \vd s \right)} f(X_t) \right] .\]
Since $b'\equiv 0$, It\^o-Tanaka formula applied to the function $B(y):=\int_0^y {b}(x) \vd x$ and the driftless $\beta$-SBM $X_t$ implies that $\PP$-a.s.~
\[
\begin{split}
	& \int_0^t  b(X_s) \vd W_s - \frac12 \int_0^t  ({b}(X_s))^2 \vd s 
	= B(X_s) -B(X_0)- \frac12 \int_0^t  ({b}(X_s))^2 + {b}'(X_s)    \vd s  - C_{b,\sigma}  L^0_t(X)
\end{split}\]
where $C_{b,\sigma}:= \frac{b_+ \sigma_- - b_- \sigma_+}{\sigma_-+\sigma_+}\in (-\infty,0)$.
Therefore, 
\[ 
\EE\!\left[ f(X_t) \right] 
= \EE_{\QR}\!\left[  \exp{\!\left( B(X_t) - B(X_0) - C_{b,\sigma} L^0_t(X) \right)} f(X_t) \right]. 
\]
Then
\[ 
\EE\!\left[ f(X_t) \right] 
\leq e^{-B(X_0)} \, \EE_{\QR}\!\left[ e^{B(X_t)} e^{ - C_{b,\sigma} L^0_t(X)}  f(X_t) \right] 
\leq e^{-B(X_0)} \, \EE_{\QR}\!\left[ e^{B(X_t)} f(X_t) \EE_{\QR}\!\left[ e^{- C_{b,\sigma} L^0_t(X)} | X_t \right] \right].
\]

To complete the proof we need to show that 
\begin{equation} \label{eq:density:key}
 \EE_{\QR}\!\left[ e^{- C_{b,\sigma} L^0_t(X)} | X_t \right] \leq 1+ \sqrt{2 \pi t} \, | C_{b,\sigma} |.
\end{equation}
Let us first aussume~\eqref{eq:density:key}. Then, for every $f$ measurable and bounded
\[ 
\EE\!\left[ f(X_t) \right] 
\leq e^{-B(X_0)} \left(1+ \sqrt{2 \pi t} \, \frac{|b_+\sigma_- - b_- \sigma_+|}{\sigma_-+\sigma_+} \right) \, \EE_{\QR}\!\left[ f(X_t) \right]
\]
and we complete the proof as in~\cite[Corollary~5]{downes2009bounds}: to consider the transition density we take derivative of the partition function of 
$X_t$. More precisely, we consider $f=\ind{(y-\varepsilon, y+\varepsilon]}$ for $\varepsilon >0$, divide for $\varepsilon $ and take the limit on both sides as $\varepsilon \to 0$. Thus,
\[ p_{(b_+,b_-)}(t,X_0,y) \leq e^{B(y)-B(x_0)} \left(1+\sqrt{ 2\pi t} \frac{|b_+\sigma_- - b_- \sigma_+|}{\sigma_-+\sigma_+}\right) p_{(0,0)}(t,X_0,y) \]
and the transition density of a $\beta$-SBM, is known (see \cite{walsh}):
\begin{equation} \label{eq:skew:density}
    p_{(0,0)}(t,x,y)
    =\frac{1}{\sqrt{2\pi t}}\exp\left(-\frac{(x-y)^2}{2t}\right)
    +\beta \sgn(y)\frac{1}{\sqrt{2\pi t}}\exp\left(-\frac{(|{x}|+|{y}|)^2}{2t}\right)
\end{equation}
which is bounded from above by $\frac{1+|\beta|}{\sqrt{2\pi t}} \exp\left(-\frac{(x-y)^2}{2t}\right)$.

To conclude, we now prove~\eqref{eq:density:key}. In what follows let $C:=|C_{b,\sigma}|$.\\
Let $\rho(t,x,y)$ be the joint transition density of the skew BM $Y_t=x + W_t + \frac{\sigma_--\sigma_+}{\sigma_-+\sigma_+} L_t^0(Y)$, $t\in [0,\infty)$, and its local time: $(Y_t,L^0_t(Y))$ (see~\cite[Proposition~1]{EMloc} or~\cite{abtww}):
\[
	\rho_t(\ell,X_0, X_t) (1-\delta_0(\ell))= 
	\frac{(1+\sgn(X_t)\frac{\sigma_--\sigma_+}{\sigma_-+\sigma_+})(\ell+ X_0 +|X_t|)}{\sqrt{2\pi t^3}} e^{-\frac{(\ell + X_0 +|X_t|)^2}{2 t}}
\]
Recall that
\[ \EE_{\QR}\!\left[ e^{C L^0_t(X)} | X_t \right] =
  \int_0^\infty e^{C \ell} \frac{\rho_t(\ell,X_0, X_t)}{p_{(0,0)}(t,X_0,X_t)} \vd \ell
\]
where $p_{(0,0)}(t,X_0, X_t)$ is given by~\eqref{eq:skew:density}.
It holds for all $c\in (0,\infty)$ that
\[
\begin{split}
&
	e^{\frac{c \ell}{t} } \frac{\rho_t(\ell,X_0, X_t)}{p_{(0,0)}(t,X_0, X_t)} (1-\delta_0(\ell)) 
\\
& 
= 
\frac{(\ell+ X_0 +|X_t|-c)+c}{t} e^{-\frac{(\ell + X_0 +|X_t| -c )^2}{2 t}}
	\frac{ (1+\sgn(X_t)\frac{\sigma_--\sigma_+}{\sigma_-+\sigma_+} )  e^{\frac{ ( X_0 +|X_t| -c)^2 - 4  X_t X_0 \ind{\{X_t>0\}}}{2 t}}}{ 1 + \sgn(X_t) \frac{\sigma_--\sigma_+}{\sigma_-+\sigma_+} e^{-2\frac{X_t X_0}{t} \ind{\{X_0 X_t>0\}}}}.
\end{split}
\]
Note that
\[
\begin{split}
&\int_0^\infty \frac{(\ell+ c_2)+c}{t} e^{-\frac{(\ell + c_2 )^2}{2 t}}  \vd \ell
= \int_{c_2}^\infty \frac{\ell+c}{t} e^{-\frac{\ell^2}{2 t}}  \vd \ell
= e^{-\frac{(c_2)^2}{2 t}} + \frac{c\sqrt{2 \pi}}{\sqrt{t}} \Phi\!\left(-\frac{c_2}{\sqrt{t}}\right)
\end{split}
\]
with $c=C t$ and $c_2 = X_0+|X_t| -C t$.
Therefore
\[
\EE_{\QR}\!\left[ e^{C L^0_t(X)} | X_t \right] 
	= e^{-\frac{(X_0+|X_t| -C t)^2}{2 t}} + C \sqrt{t} \sqrt{2 \pi} \Phi\!\left(-\frac{X_0+|X_t| -C t}{\sqrt{t}}\right) \leq 1+ \sqrt{2\pi t} \, C.
\]
This completes the proof.
\end{proof}

\paragraph*{Acknowledgements.} We are grateful to A.~Lejay and R.~Peyre for discussions.
PP acknowledges financial support from University of Rome Tor Vergata via Project AIF - E83C22002610005.

\addcontentsline{toc}{section}{Bibliography}
\bibliographystyle{abbrvnat}
\bibliography{biblio}

\begin{thebibliography}{50}
\providecommand{\natexlab}[1]{#1}
\providecommand{\url}[1]{\texttt{#1}}
\expandafter\ifx\csname urlstyle\endcsname\relax
  \providecommand{\doi}[1]{doi: #1}\else
  \providecommand{\doi}{doi: \begingroup \urlstyle{rm}\Url}\fi

\bibitem[Amorino and Gloter(2020)]{amorino_glotier}
C.~Amorino and A.~Gloter.
\newblock Contrast function estimation for the drift parameter of ergodic jump
  diffusion process.
\newblock \emph{Scandinavian Journal of Statistics}, 47\penalty0 (2):\penalty0
  279--346, 2020.

\bibitem[Ang and Timmermann(2012)]{ang}
A.~Ang and A.~Timmermann.
\newblock {Regime Changes and Financial Markets}.
\newblock \emph{Annual Review of Financial Economics}, 4:\penalty0 313--337,
  2012.

\bibitem[Appuhamillage et~al.(2011)Appuhamillage, Bokil, Thomann, Waymire, and
  Wood]{abtww}
T.~Appuhamillage, V.~Bokil, E.~Thomann, E.~Waymire, and B.~Wood.
\newblock Occupation and local times for skew {B}rownian motion with
  applications to dispersion across an interface.
\newblock \emph{Ann. Appl. Probab.}, 21\penalty0 (1):\penalty0 183--214, 2011.

\bibitem[Bass and Chen(2005)]{bass2005one}
R.~F. Bass and Z.-Q. Chen.
\newblock One-dimensional stochastic differential equations with singular and
  degenerate coefficients.
\newblock \emph{Sankhy{\=a}: The Indian Journal of Statistics}, pages 19--45,
  2005.

\bibitem[Ben~Alaya and Kebaier(2013)]{alaya_kebaier_2013}
M.~Ben~Alaya and A.~Kebaier.
\newblock {Asymptotic Behavior of the Maximum Likelihood Estimator for Ergodic
  and Nonergodic Square-Root Diffusions}.
\newblock \emph{Stochastic Analysis and Applications}, 31\penalty0
  (4):\penalty0 552--573, 2013.

\bibitem[Borodin(1986)]{Borodin}
A.~N. Borodin.
\newblock On the character of convergence to brownian local time. ii.
\newblock \emph{Probability Theory and Related Fields}, 72:\penalty0 251--277,
  1986.

\bibitem[Borodin and Salminen(2015)]{borodin2015handbook}
A.~N. Borodin and P.~Salminen.
\newblock \emph{Handbook of Brownian motion-facts and formulae}.
\newblock Springer Science \& Business Media, 2015.

\bibitem[Brooks et~al.(2011)Brooks, Gelman, Jones, and
  Meng]{brooks2011handbook}
S.~Brooks, A.~Gelman, G.~Jones, and X.-L. Meng.
\newblock \emph{Handbook of Markov Chain Monte Carlo}.
\newblock CRC press, 2011.

\bibitem[{Buckner} et~al.(2024){Buckner}, {Dowd}, and {Hulley}]{buckner}
D.~{Buckner}, K.~{Dowd}, and H.~{Hulley}.
\newblock {Arbitrage Problems with Reflected Geometric Brownian Motion}.
\newblock \emph{Finance Stoch}, 28, Jan. 2024.

\bibitem[Chen et~al.(2011)Chen, So, and Liu]{Chen:2011bk}
C.~W.~S. Chen, M.~K.~P. So, and F.-C. Liu.
\newblock A review of threshold time series models in finance.
\newblock \emph{Statistics and its Interface}, 4\penalty0 (2):\penalty0
  167--181, 2011.

\bibitem[Christensen and Strauch(2023)]{ChristensenStrauchAAP}
S.~Christensen and C.~Strauch.
\newblock {Nonparametric learning for impulse control problems -- Exploration
  vs. exploitation}.
\newblock \emph{The Annals of Applied Probability}, 33\penalty0 (2):\penalty0
  1569 -- 1587, 2023.

\bibitem[Christensen et~al.(2023)Christensen, Strauch, and
  Trottner]{ChristensenStrauchTrottnerBernoulli}
S.~Christensen, C.~Strauch, and L.~Trottner.
\newblock Learning to reflect: A unifying approach for data-driven stochastic
  control strategies.
\newblock \emph{Bernoulli}, 2023.

\bibitem[Crimaldi and Pratelli(2005)]{crimaldi}
I.~Crimaldi and L.~Pratelli.
\newblock Convergence results for multivariate martingales.
\newblock \emph{Stochastic Process. Appl.}, 115\penalty0 (4):\penalty0
  571--577, 2005.

\bibitem[Decamps et~al.(2006)Decamps, Goovaerts, and Schoutens]{interestrate}
M.~Decamps, M.~Goovaerts, and W.~Schoutens.
\newblock Self exciting threshold interest rates models.
\newblock \emph{Int. J. Theor. Appl. Finance}, 9\penalty0 (7):\penalty0
  1093--1122, 2006.

\bibitem[Dieker and Gao(2013)]{dieker2013}
A.~B. Dieker and X.~Gao.
\newblock {Positive recurrence of piecewise Ornstein-Uhlenbeck processes and
  common quadratic Lyapunov functions}.
\newblock \emph{Ann. Appl. Probab.}, 23\penalty0 (4):\penalty0 1291--1317, 08
  2013.

\bibitem[Dong and Wong(2017)]{DongWong}
F.~Dong and H.~Y. Wong.
\newblock {Variance Swaps under the Threshold Ornstein-Uhlenbeck Model}.
\newblock \emph{Appl. Stoch. Model. Bus. Ind.}, 33\penalty0 (5):\penalty0
  507--521, Sept. 2017.

\bibitem[Downes(2009)]{downes2009bounds}
A.~N. Downes.
\newblock Bounds for the transition density of time-homogeneous diffusion
  processes.
\newblock \emph{Statistics \& probability letters}, 79\penalty0 (6):\penalty0
  835--841, 2009.

\bibitem[{\'E}tor{\'e} and Martinez(2013)]{EMloc}
P.~{\'E}tor{\'e} and M.~Martinez.
\newblock Exact simulation of one-dimensional stochastic differential equations
  involving the local time at zero of the unknown process.
\newblock \emph{Monte Carlo Methods Appl.}, 19\penalty0 (1):\penalty0 41--71,
  2013.

\bibitem[Florens-Zmirou(1993)]{florens-zmirou}
D.~Florens-Zmirou.
\newblock On estimating the diffusion coefficient from discrete observations.
\newblock \emph{J. Appl. Probab.}, 30\penalty0 (4):\penalty0 790--804, 1993.

\bibitem[Gairat and Shcherbakov(2016)]{gairat}
A.~Gairat and V.~Shcherbakov.
\newblock Density of skew {B}rownian motion and its functionals with
  application in finance.
\newblock \emph{Mathematical Finance}, 26\penalty0 (4):\penalty0 1069--1088,
  2016.

\bibitem[Giorgi et~al.(2023)Giorgi, Herzel, and Pigato]{ghp2023}
F.~Giorgi, S.~Herzel, and P.~Pigato.
\newblock A reinforcement learning algorithm for trading commodities.
\newblock \emph{Applied Stochastic Models in Business and Industry},
  n/a\penalty0 (n/a), 2023.

\bibitem[Hottovy and Stechmann(2015)]{Hottovy}
S.~Hottovy and S.~N. Stechmann.
\newblock Threshold models for rainfall and convection: Deterministic versus
  stochastic triggers.
\newblock \emph{SIAM Journal on Applied Mathematics}, 75\penalty0 (2):\penalty0
  861--884, 2015.

\bibitem[Hu and Xi(2022)]{HuXi}
Y.~Hu and Y.~Xi.
\newblock Parameter estimation for threshold ornstein--uhlenbeck processes from
  discrete observations.
\newblock \emph{Journal of Computational and Applied Mathematics},
  411:\penalty0 114264, 2022.

\bibitem[It\^o and McKean(1996)]{ito}
K.~It\^o and H.~McKean.
\newblock \emph{Diffusion Processes and their Sample Paths: Reprint of the 1974
  Edition}.
\newblock Classics in Mathematics. Springer Berlin Heidelberg, 1996.

\bibitem[Jacod(1998)]{j1}
J.~Jacod.
\newblock Rates of convergence to the local time of a diffusion.
\newblock \emph{Ann. Inst. H. Poincar\'e Probab. Statist.}, 34\penalty0
  (4):\penalty0 505--544, 1998.

\bibitem[Jacod and Protter(2012)]{jp}
J.~Jacod and P.~Protter.
\newblock \emph{Discretization of processes}, volume~67 of \emph{Stochastic
  Modelling and Applied Probability}.
\newblock Springer, Heidelberg, 2012.

\bibitem[Jacod and Shiryaev(2003)]{js}
J.~Jacod and A.~N. Shiryaev.
\newblock \emph{Limit theorems for stochastic processes}, volume 288 of
  \emph{Grundlehren der Mathematischen Wissenschaften}.
\newblock Springer-Verlag, Berlin, second edition, 2003.

\bibitem[Keilson and Wellner(1978)]{kw}
J.~Keilson and J.~A. Wellner.
\newblock Oscillating {B}rownian motion.
\newblock \emph{J. Appl. Probability}, 15\penalty0 (2):\penalty0 300--310,
  1978.

\bibitem[Kessler(1997)]{kessler}
M.~Kessler.
\newblock {Estimation of an Ergodic Diffusion from Discrete Observations}.
\newblock \emph{Scandinavian Journal of Statistics}, 24\penalty0 (2):\penalty0
  211--229, 1997.

\bibitem[Kutoyants(2012)]{kutoyants2012}
Y.~A. Kutoyants.
\newblock On identification of the threshold diffusion processes.
\newblock \emph{Ann. Inst. Statist. Math.}, 64\penalty0 (2):\penalty0 383--413,
  2012.

\bibitem[Le~Gall(1985)]{legall}
J.-F. Le~Gall.
\newblock One-dimensional stochastic differential equations involving the local
  times of the unknown process.
\newblock \emph{Stochastic Analysis. Lecture Notes Math.}, 1095:\penalty0
  51--82, 1985.

\bibitem[Lejay and Pigato(2018)]{LP}
A.~Lejay and P.~Pigato.
\newblock Statistical estimation of the {O}scillating {B}rownian {M}otion.
\newblock \emph{Bernoulli}, 24\penalty0 (4B):\penalty0 3568--3602, 2018.

\bibitem[Lejay and Pigato(2019)]{lp1}
A.~Lejay and P.~Pigato.
\newblock A threshold model for local volatility: evidence of leverage and mean
  reversion effects on historical data.
\newblock \emph{International Journal of Theoretical and Applied Finance},
  22\penalty0 (4):\penalty0 1950017, 2019.

\bibitem[Lejay and Pigato(2020)]{lp2}
A.~Lejay and P.~Pigato.
\newblock Maximum likelihood drift estimation for a threshold diffusion.
\newblock \emph{Scandinavian Journal of Statistics}, 47\penalty0 (3):\penalty0
  609--637, 2020.

\bibitem[Lejay et~al.(2019)Lejay, Mordecki, and Torres]{lmt1}
A.~Lejay, E.~Mordecki, and S.~Torres.
\newblock Two consistent estimators for the skew {B}rownian motion.
\newblock \emph{ESAIM Probab. Stat.}, 23:\penalty0 567--583, 2019.

\bibitem[Lipton(2018)]{lipton:2018}
A.~Lipton.
\newblock {Oscillating Bachelier and Black-Scholes Formulas}.
\newblock In \emph{Financial Engineering}. World Scientific, 2018.

\bibitem[Lipton and Sepp(2011-10)]{liptonsepp}
A.~Lipton and A.~Sepp.
\newblock Filling the gaps.
\newblock \emph{Risk Magazine}, pages 66--71, 2011-10.

\bibitem[Mazzonetto(2019)]{mazzonetto2019rates}
S.~Mazzonetto.
\newblock Rates of convergence to the local time of oscillating and skew
  brownian motions.
\newblock \emph{arXiv preprint arXiv:1912.04858}, 2019.

\bibitem[Mazzonetto and Pigato(2024)]{mazzonetto2020drift}
S.~Mazzonetto and P.~Pigato.
\newblock Drift estimation of the threshold ornstein-uhlenbeck process from
  continuous and discrete observations.
\newblock \emph{Statistica Sinica}, 34:\penalty0 313--336, 2024.

\bibitem[Meyn and Tweedie(2009)]{MeynTweedie}
S.~Meyn and R.~L. Tweedie.
\newblock \emph{Markov Chains and Stochastic Stability}.
\newblock Cambridge University Press, USA, 2nd edition, 2009.
\newblock ISBN 0521731828.

\bibitem[Mota and Esqu{\'{\i}}vel(2014)]{motaesquivel}
P.~P. Mota and M.~L. Esqu{\'{\i}}vel.
\newblock On a continuous time stock price model with regime switching, delay,
  and threshold.
\newblock \emph{Quant. Finance}, 14\penalty0 (8):\penalty0 1479--1488, 2014.

\bibitem[Pai and Pedersen(1999)]{pai}
J.~Pai and H.~Pedersen.
\newblock {Threshold Models of the Term Structure of Interest Rate}.
\newblock In \emph{Joint day Proceedings Volume of the XXXth International
  ASTIN Colloquium/9th International AFIR Colloquium, Tokyo, Japan}, pages
  387--400. 1999.

\bibitem[Pigato(2019)]{pigato}
P.~Pigato.
\newblock { Extreme at-the-money skew in a local volatility model}.
\newblock \emph{Finance and Stochastics}, 23:\penalty0 827--859, 2019.

\bibitem[Portenko(1994)]{port1}
N.~I. Portenko.
\newblock The development of {I}. {I}. {G}ikhman's idea concerning the methods
  for investigating local behavior of diffusion processes and their weakly
  convergent sequences.
\newblock \emph{Teor. \u{I}mov\={\i}r. Mat. Stat.}, \penalty0 (50):\penalty0
  6--21, 1994.

\bibitem[R{\'e}nyi(1963)]{renyi}
A.~R{\'e}nyi.
\newblock On stable sequences of events.
\newblock \emph{Sankhy\=a Ser. A}, 25:\penalty0 293--302, 1963.

\bibitem[Su and Chan(2015)]{su2015}
F.~Su and K.-S. Chan.
\newblock Quasi-likelihood estimation of a threshold diffusion process.
\newblock \emph{J. Econometrics}, 189\penalty0 (2):\penalty0 473--484, 2015.

\bibitem[Su and Chan(2017)]{su2017}
F.~Su and K.-S. Chan.
\newblock Testing for threshold diffusion.
\newblock \emph{J. Bus. Econom. Statist.}, 35\penalty0 (2):\penalty0 218--227,
  2017.

\bibitem[Tong(2011)]{Tong:2011ud}
H.~Tong.
\newblock Threshold models in time series analysis --- 30 years on.
\newblock \emph{Statistics and its Interface}, 4, 2011.

\bibitem[Walsh(1978)]{walsh}
J.~B. Walsh.
\newblock A diffusion with discontinuous local time.
\newblock In \emph{Temps locaux}, volume 52-53, pages 37--45. Soci{\'e}t{\'e}
  Math{\'e}matique de France, 1978.

\bibitem[Yu et~al.(2020)Yu, Tsai, and Rachinger]{YU2020}
T.-H. Yu, H.~Tsai, and H.~Rachinger.
\newblock Approximate maximum likelihood estimation of a threshold diffusion
  process.
\newblock \emph{Computational Statistics \& Data Analysis}, 142:\penalty0
  106823, 2020.

\end{thebibliography}

\end{document}